\documentclass[final]{siamltex}
\usepackage{amsmath, amsfonts,color}
\usepackage[english]{babel}
\usepackage{dsfont}
\usepackage{amssymb}
\usepackage{mathabx}
\usepackage{enumitem}
\usepackage{graphicx,epsfig}
\usepackage{tikz}
\usepackage{float}
\usepackage{color}
\usepackage{capt-of}
\usepackage{tikz}
\usepackage{hyperref}

\newcommand{\R}{{\mathbb R}}

\newcommand{\EE}{{\mathbb E}}
\newcommand{\one}{{\mathds{1}}}
\newcommand{\PP}{{\mathbb P}}

\newcommand{\cS}{{\mathcal S}}
\newcommand{\cP}{{\mathcal P}}
\newcommand{\cN}{{\mathcal N}}

\newcommand{\cF}{{\mathcal F}}
\newcommand{\cV}{{\mathcal V}}

\newcommand{\dt}{{\Delta t}}
\def\ds{\displaystyle}
 \definecolor{mypurple}{RGB}{140,0,255}
\definecolor{myred}{rgb}{255,0,0}

\definecolor{mydarkturquoise}{RGB}{0,206,209}
\definecolor{mydeeppink}{RGB}{255,20,147}
\definecolor{darkblue}{RGB}{0,0,140}
\definecolor{blue2}{RGB}{0,0,0}
\definecolor{middleblue}{RGB}{0,0,71}
\definecolor{light-gray}{gray}{0.9}

\definecolor{ProcessBlue}{cmyk}{1,0,0,0.40}
\definecolor{Black}{cmyk}{0,0,0,1}
\definecolor{Red}{cmyk}{0,1,1,0.2}
\definecolor{Green}{cmyk}{0.9,0,1,0}
\definecolor{Orange}{cmyk}{0,0.61,0.87,0.5}
\definecolor{Fuchsia}{cmyk}{0.47,0.91,0,0.06}
\definecolor{PineGreen}{cmyk}{0.92,0,0.59,0.30}

\numberwithin{equation}{section}
 \def\theequation{\thesection.\arabic{equation}}
 \newtheorem{remark}{Remark}[section]

\begin{document}

\title{Simulating numerically the Krusell-Smith model with neural networks}
\author{Yves Achdou \thanks { Universit{\'e}  Paris-Cit{\'e} and  Sorbonne Universit{\'e}, CNRS, Laboratoire Jacques-Louis Lions, 
(LJLL), F-75006 Paris, France, achdou@ljll-univ-paris-diderot.fr} \and
Jean-Michel Lasry \thanks{Universit{\'e} Paris-Dauphine, France} \and
Pierre-Louis Lions\thanks{Coll{\`e}ge de France,  Paris, France } 
}

\maketitle
\begin{abstract}
  The celebrated Krusel-Smith growth model is an important example of a Mean Field Game with a common noise.
  The Mean Field Game is encoded in the master equation, a partial differential equation satisfied by the value of the game which depends on the whole  distribution of states. The latter equation is therefore  posed in an infinite dimensional space. This makes the numerical simulations quite challenging. However, Krusell and Smith conjectured that the
  value function of the game mostly depends on the state distribution through low dimensional quantities.   In this paper,  we wish to propose a numerical method for approximating the solutions of the master equation arising in Krusell-Smith model, and for adaptively identifying  low-dimensional variables which retain
an important part of the information. This new numerical framework is based on a semi-Lagrangian method and uses neural networks as an important ingredient.
\end{abstract}
\section{Introduction}

Krusell-Smith model, see \cite{krusell-smith}, is a celebrated growth model in macroeconomics. In this model, the agents are households whose wealth and productivity are heterogeneous; they aim at maximizing some criteria involving their consumption and are bound by a borrowing limit. Krusell-Smith model differs from  the previous Aiyagari-Bewley-Huggett models, see \cite{huggett, aiyagari,bewley}, in which the productivities of the agents are subject to some idiosynchratic noise,  because it incorporates random shocks which affect the whole economy.  As we shall see below, the introduction of a common noise implies a major additional difficulty and  a real challenge in macroeconomics: roughly speaking,  the optimal value and the optimal strategy of a single agent cannot be simply expressed in terms of his own wealth and productivity, but also depends on the distribution of states of  other agents, which is a quantity in an infinite dimensional space.

However, in \cite{krusell-smith}, Krusell and Smith made the important conjecture that the optimal value  depends on the distribution of states mostly through a finite dimensional information, and even through a single number which is besides a nonlinear function of the latter distribution. This would mean that the optimal strategy of the agents mostly depends on a single common parameter.

The model of Krusell and Smith and the related open questions have been formulated in the language of macro-economics and were lacking of a precise mathematical formulation.

Later and independently, two authors of the present paper  proposed the mathematical theory of Mean Field Games  ({\sl MFGs} in short), see ~\cite{PLL-CDF,MR2269875,MR2271747,MR2295621}, which aims at modelling dynamical equilibria for large populations, as it is the case, for instance, when  studying deterministic or stochastic differential
games (Nash equilibria) as the number of agents tends to infinity. 
In that case, one assumes that the rational agents are  indistinguishable and  individually have a negligible influence on the game, 
and that each individual strategy is influenced for example by some averages of quantities depending on the states (or the controls) of the other agents. In 2009, the last two authors interacted with R. Lucas who drew their attention to Krusell-Smith model and the related mathematical challenges. Although they already knew how to write, in particular for Mean Field Games with common noise, what is known nowadays as  the {\sl master equation}, the Krusell-Smith model was crucial in their placing this equation at the center of the theory. Indeed, the Krusell-Smith model can be seen as a Mean Field Game in which the agents interact through aggregate quantities, and that naturally leads to such a master equation, see for instance \cite{MR3268061}. Note that the terminology {\sl master equation} is inspired from statistical physics. It is a partial derivative equation ({\sl PDE} in short) satisfied by the optimal value  function of the considered  Nash equilibrium with a continuum of agents: since, as in Krusell-Smith model, the latter optimal value  depends on the whole distribution of states, the  PDE is posed in an infinite dimensional space, and new mathematical notions  are needed to give a meaning to the derivatives with respect to the probability measure associated with the distribution of states, see \cite{PLL-CDF,MR3967062}.

The master equation has been useful to give a precise mathematical meaning to Krusell-Smith model and conjecture, but has not been used yet to check whether the latter is true; the reason for that is the infinite dimensionality of the PDE, which makes it very difficult to apply numerical methods. Even though, since 2009, many progress have been made in the mathematical anlysis of the master equation, and in using  the latter for modeling in economics and other social sciences, the problem of finding efficient numerical approximations  remains open.

The purpose of the present paper is to propose a numerical method for approximating the solutions of the master equation arising in Krusell-Smith model and apply it to address the above mentioned conjecture.

This new numerical framework is based on a semi-Lagrangian method (see the appendix by M. Falcone in the book by M. Bardi and I. Capuzzo-Dolcetta, \cite{MR1484411}, for semi-Lagrangian methods in the context of
optimal control theory) and uses as an important ingredient neural networks to cope with high dimensionality, but also some mathematical understanding of its solution, see \cite{MR4365976}. Besides, let us note that the fixed point formulation  arising from the abovementioned numerical method might also appear easier to understand for readers who are not familiar with infinite dimensional PDEs.

 \noindent {\bf Acknowledgments.}
This research  was partially supported by the chair Finance and Sustainable Development and FiME Lab (Institut Europlace de Finance).

\section{The Krusell-Smith model}

We consider households (named agents hereafter) which are heterogeneous in  their wealth (or capital) $x$ and productivity $y$. The dynamics of the wealth of a given agent is given by 
\begin{equation*}     
 dx_t = (r_t x_t+ w_t y_t  -c_t ) dt,
    \end{equation*}
where 
\begin{itemize}
\item $y_t$ is the productivity of the agent (a state variable)
\item $c_t$ is the consumption (the control variable)
\item $r_t$ is the interest rate (common to all agents)
\item $w_t$ is the unitary salary (common to all agents)
\end{itemize}
 The productivity $y_t$ is a two-state Poisson process with intensities $\lambda_1$ and $\lambda_2$, i.e.
$y_t \in \{ y_1, y_2\}$ with $y_1<y_2$, and 
      \begin{eqnarray*} 
          &\PP\left(y_{t+\dt} =y_1\right| \left.  y_{t} =y_1 \right)= 1- \lambda_1 \dt +o(\dt),\\
          &\PP\left(y_{t+\dt} =y_2\right| \left.  y_{t} =y_1 \right)= \lambda_1 \dt +o(\dt) ,
 \\
          &   \PP\left(y_{t+\dt} =y_2\right| \left.  y_{t} =y_2 \right)= 1-\lambda_2 \dt +o(\dt), \\  
          & \PP\left(y_{t+\dt} =y_1\right| \left.  y_{t} =y_2 \right)= \lambda_2 \dt +o(\dt).
      \end{eqnarray*}
      We assume that  the  random processes describing the productivities  of the agents are all independent (idiosynchratic noises).\\
Recall that a negative wealth means debts; there is a {\sl borrowing constraint}: 
the wealth of a given household cannot be less than a given {\sl borrowing limit}  $\underline{x}$. 
In the terminology of control theory,  $  x_t\ge  \underline{x}$ is a constraint on the state variable.\\
To determine the interest rate $r_t$ and the level of wages $w_t$, we assume that the production of the economy is  described  
 by the following Cobb-Douglas law:
\[  F_t(X_t,Y_t)=  A_t X_t^\alpha Y_t^{1-\alpha}  ,\]
where
\begin{itemize}
  \item the exponent $\alpha$ lies in $(0,1)$
  \item $A_t$ is a  noisy productivity factor (the noise affects the whole economy and is independent from the productivities of the individuals)
  \item $X_t= \int_{x\ge \underline  x} \int_{y\in \{y_1,y_2\}}  x dm (t,x,y) $ is the aggregate capital
  \item $Y_t=   \int_{x\ge \underline  x} \int_{y\in \{y_1,y_2\}}  y dm (t,x,y)$ is  the  aggregate  labor
  \item The distribution $m(t,\cdot,\cdot)$ of the pairs $(x_t,y_t)$  is a probability measure on $[\underline x,+\infty)\times \{y_1, y_2\}$.
  \end{itemize}
 The level of wages $w_t$ and the interest rate $r_t$ are obtained by  the equilibrium relation
  \begin{displaymath}
    (X_t,Y_t)= {\rm{argmax}}    \Bigl( F_t(X,Y)- (r_t+\delta ) X -w_t Y \Bigr),
  \end{displaymath}
  where $\delta$ is the rate of depreciation of the capital. This implies that 
  \begin{displaymath}
  r_t=\partial _X F_t (X_t,Y_t)-\delta  = \alpha A_t \frac{Y_t^{1-\alpha} }{X_t^{1-\alpha} } -\delta , \quad\quad w_t = \partial _Y F_t (X_t,Y_t) =
 (1- \alpha )A_t \frac{X_t^{\alpha} }{Y_t^{\alpha} } .
\end{displaymath}
We assume that $A_t$ is a two-state Poisson process independent from the noises affecting the productivity of the agents,
 with intensities $\mu_1$ and $\mu_2$, i.e.
 i.e.
$A_t \in \{ A_1, A_2\}$ with $A_1<A_2$, and 
      \begin{eqnarray*} 
           &\PP\left(A_{t+\dt} =A_1\right| \left.  A_{t} =A_1 \right)= 1- \mu_1 \dt +o(\dt), \\  &  \PP\left(A_{t+\dt} =A_2\right| \left.  A_{t} =A_1 \right)= \mu_1 \dt +o(\dt) ,\\
          &   \PP\left(A_{t+\dt} =A_2\right| \left.  A_{t} =A_2 \right)= 1-\mu_2 \dt +o(\dt), \\   
&          \PP\left(A_{t+\dt} =A_1\right| \left.  A_{t} =A_2 \right)= \mu_2 \dt +o(\dt).
      \end{eqnarray*}
      In what follows, we set
\begin{equation}\label{eq:5}
  r_i(m)= \alpha A_i \frac{\left(\int_{x\ge \underline  x} \int_{y\in \{y_1,y_2\}}  y dm (x,y) \right)^{1-\alpha} }{\left(\int_{x\ge \underline  x} \int_{y\in \{y_1,y_2\}}  x dm (x,y) \right)^{1-\alpha} } -\delta,
\end{equation}
and 
\begin{equation}\label{eq:6}
  w_i(m)= (1-\alpha) A_i \frac{\left(\int_{x\ge \underline  x} \int_{y\in \{y_1,y_2\}}  x dm (x,y) \right)^{\alpha} }  {\left(\int_{x\ge \underline  x} \int_{y\in \{y_1,y_2\}}  y dm (x,y) \right)^{\alpha} }.
\end{equation}
An agent solves the optimal control problem
    \begin{equation*}
      \max_{\{c_t\}}  \mathbb{E}\int_0^\infty e^{-\rho t} u(c_t)dt \quad \hbox{subject to }
      \left\{
        \begin{array}[c]{rcl}
          dx_t &=& (w_t y_t + r_t x_t - c_t)dt,\\
          x_t  &\geq& \underline{x},    
        \end{array}
      \right.   
    \end{equation*}
    where
   \begin{itemize}   
     \item $\rho$ is a positive discount factor
     \item $u$ is a  utility function,  strictly increasing and strictly concave, e.g. the CRRA (constant relative risk aversion) utility:
       \begin{equation}
         \label{eq:1}
        u(c)=c^{1-\gamma}/(1-\gamma),\qquad \gamma>0.
       \end{equation}
     \end{itemize}
     The introduction of aggregate shocks (on $A_t$) creates a major difficulty: in contrast with the case without aggregate uncertainty, it becomes necessary to include the entire distribution of productivity and wealth $m$ as a state variable in the optimal control problem of the individuals. This distribution is now itself a random variable and hence calendar time $t$ is no longer a sufficient statistic to describe the behavior of the system.

The aggregate state is $(A_i,m),i=1,2$ and the individual state is $(x,y)$. The value of an individual agent 
when  $A_t=A_i$, $i=1,2$, is $v_i(x,y,m)$.

The master equations satisfied by the value functions $v_i$ are posed in $(\underline x, +\infty)\times \{y_1,y_2\} \times \PP( [\underline x, +\infty)\times \{y_1,y_2\})$ and read as follows: for $i=1,2$, $\bar \imath = 3 -i $, $j=1,2$, $\bar \jmath = 3 -j $, 
\begin{equation*}
  \begin{array}[c]{l}
    \begin{aligned}[t]
0 = &    \lambda_j  (v_i(x,y_{\bar \jmath})-v_i(x,y_j) )   + (w_i(m)y_j + r_i(m) x)\partial_x v_i(x,y_j)  +H(\partial_x v_i (x,y_j))\\
&  - \rho v_i(x,y_j)\ + \mu_i(v_{\bar \imath}(x,y_j) - v_i(x,y_j))\\ 
&+  \sum_{\ell=1}^2 \int_{\hat x}    T[m,\partial_x v_{i}](\hat x,y_\ell) \frac{\delta v_i}{\delta m}(x,y_j, \hat x, y_\ell)  d\hat x ,
\end{aligned}
    \\
\hbox{with} \\
  \begin{aligned}[t]
T[m,\partial_x v_{i}] (\hat x, y_\ell) = & \lambda_{\bar \ell} m(\hat x,y_{\bar \ell})-\lambda_\ell m(\hat x,y_\ell)  
-\partial_x \Bigl(   (w_i(m)y_\ell + r_i(m) \cdot ) m (\cdot, y_\ell)\Bigr) (\hat x)\\
&- \partial_x \Bigl(\partial_p H(\partial_x v_i (\cdot, y_\ell))  m(\cdot,y_\ell) \Bigr)(\hat x),
\end{aligned}
  \end{array}
\end{equation*}
and
\begin{equation}\label{eq:H}
H(p)= \max_{c \ge 0}    \left( -pc+u(c) \right).
\end{equation}
With $u$  given by~(\ref{eq:1}),
\begin{equation*}
  H(p)=\left\{
      \begin{aligned}[t]
       \frac \gamma {1-\gamma} p^{1-\frac 1 \gamma},\quad \quad & \hbox{if }\quad p>0,\\
       +\infty, \quad \quad  & \hbox{if }\quad p\le 0.
     \end{aligned}
   \right.
\end{equation*}

It may be more convenient to describe the value function by four functions on $[\underline x, +\infty)\times \PP( [\underline x, +\infty) \times \{y_1,y_2\}   )$, namely $v_{i,j}$, $i,j=1,2$. The distribution of states 
is  then given by two measures on $[\underline x, +\infty)$, namely $m_{j}$, $j=1,2$. The master equation then takes the form:
\begin{eqnarray}     \label{eq:2}
0 &= &   \ds  \lambda_j  (v_{i,\bar \jmath} -v_{i,j} )   + (w_i(m)y_j + r_i(m) x)\partial_x v_{i,j}  +H(\partial_x v_{i,j})\\
 \nonumber && \ds  - \rho v_{i,j} + \mu_i(v_{\bar \imath,j} - v_{i,j})
+  \sum_{\ell=1}^2 \int_{\hat x}    T_{i,\ell}[m,\partial_x v_{i,\ell}](\hat x) \frac{\delta v_{i,j}}{\delta m_\ell}(\cdot, \hat x)  d\hat x ,
\\ \nonumber  \quad \quad\quad &&
\\\label{eq:3}
 T_{i,\ell}[m,\partial_x v_{i,\ell}] (\hat x) &= & \ds \lambda_{\bar \ell} m_{\bar\ell}(\hat x)-\lambda_\ell m_{\ell}(\hat x)  
-\partial_x \Bigl(   (w_i(m)y_\ell + r_i(m) \cdot ) m_{\ell} \Bigr) (\hat x)\\ \nonumber
&&\ds - \partial_x \Bigl(\partial_p H(\partial_x v_{i,\ell} )  m_\ell \Bigr)(\hat x).
\end{eqnarray}

\section{The semi-Lagrangian method}
\label{sec:semi-lagr-meth}
\subsection{The dynamic programming principle applied to the Nash MFG-equilibrium}\label{sec:dynam-progr-princ}
The dynamic programming principle applied to the mean field Nash equilibrium implies that for a small  time step $\Delta t$,
\begin{equation}
  \label{eq:dyn_prog}
    \begin{aligned}[t]
&  v_{i,j}(x,m)\approx \\
&\EE_{i,j}\left(
  \begin{array}[c]{l}
\Delta t  u \left(c^*_{i,j}(x,m)\right) +\\
  e^{-\rho \Delta t}   v_{i(\Delta t),j(\Delta t)}\Bigl( x+ \Delta t (w_i(m)y_j + r_i(m) x -c^*_{i,j}(x,m))  , m^*(\Delta t)    \Bigr)
  \end{array}
\right)    
\end{aligned}
\end{equation}
with an error of the order of $\Delta t$, where
$\EE_{ij} (X)$ stands to the probability of $X$ conditionned to  $A(t=0)=A_i$ and  $y(t=0)=y_j$, and 
\begin{itemize}
\item $c^*_{i,j}(x,m) $ is the optimal consumption, associated to the optimal savings policy $s^*_{i,j}(x,m) = w_i(m)y_j + r_i(m) x -c^*_{i,j}(x,m)$.
To take the state constraint into account, we set 
\begin{equation}\label{eq:7}
  c^*_{i,j}(x,m) =\min\left( -H_p\left( \partial_x  v_{i,j}(x,m) \right) ,    \frac {x- \underline x} {\Delta t  }  +w_i(m)y_j + r_i(m) x  \right)
\end{equation}
\item The probability measure $m^*(\Delta t )$ is the distribution of  $(x+ \Delta t s^*_{i,j}(x,m), y(\Delta t))$,
where the optimal policy $s^*_{i,j}(x,m)$ is defined above. 
\end{itemize}

\begin{remark}
  The dynamic programming equation  \eqref{eq:dyn_prog} might be easier to understand than the master equation for readers who are not familiar with infinite dimensional PDEs.
\end{remark}

Taylor expansions lead to the following discrete version of the master equation (\ref{eq:2})-(\ref{eq:3})
\begin{equation}\label{eq:4}
   \begin{aligned}[t]
    0=& (1+\rho \Delta t ) v_{i,j}(x,m) - v_{i,j}\Bigl(x+ \Delta t (w_i(m)y_j + r_i(m) x -c^*_{i,j}(x,m))  , m^*(\Delta t)    \Bigr)\\
&+\lambda_j \Delta t \left(  v_{i,j}(x,m)  -  v_{i,\bar \jmath}(x,m) \right) + \mu_i \Delta t \left(  v_{i,j}(x,m)  -  v_{\bar  \imath,j}(x,m) \right)\\
&+\Delta t u \left(c^*_{i,j}(x,m)\right).
\end{aligned}
\end{equation}
An important difficulty lies in the approximation of $m^*(\Delta t)$.
\subsection{Approximation of $m$ and its transported version}\label{sec:approximation-m-its}

\subsubsection{Approximation of $m$}
For computational purposes, we restrict ourselves to probability  measures $m$
belonging to a finite dimensional space; let $d$ be the  dimension of this space.
The approximation of the value function $v_{i,j}$ is therefore a function defined on $[\underline x, +\infty)\times \R^d$.\\
It is now well known that with Aiyagari and Krusell-Smith models, the distribution of capital may have a Dirac mass  
at the borrowing limit $\underline x$, see \cite{MR3268061,MR4365976}. Our space of discrete measures  therefore contain  Dirac masses at $(\underline x, y_j)$, $j=1,2$. The interpretation of these Dirac masses having
positive coefficients  is that the credit constraint is biding for a non zero percentage of the agents.
\\
We  also artificially truncate the support of the measure to the  bounded interval $[\underline x,\overline x]$ where  $\overline x$ is a sufficiently large 
positive number. To avoid loosing any mass after the measure is transported,  our space of discrete measures also contain Dirac masses at  
$(\overline x, y_j)$, $j=1,2$.
\\
Except for these four Dirac masses, our discrete measures have piecewise constant densities on $(\underline x, \overline x)$, (i.e. the latter interval 
is partitioned  into subintervals in which the density is constant, and the chosen partition depends on $j$). 
\\
Consider two subdivisions of $[\underline x, \overline x]$ associated respectively to the two increasing families: $(x_{k,j})_{k=0,\dots , K_j}$, $ j=1,2$:
\begin{displaymath}
  \underline x = x_{0,j}< x_{1,j}< \dots< x_{K_j,j}= \overline x   ,
\end{displaymath}
and set  $h_{k,j}= x_{k+1,j}-x_{k,j}$, $j=0, \dots, K_j-1$. 
The discrete  probability measures $m$ on $[\underline x, +\infty) \times \{y_1,y_2\}$ are of the form 
 $m=\sum_{j=1}^2 m_j \otimes \delta_{y=y_j} $, where
\begin{displaymath}
  m_j= \alpha_{-1,j}\delta_{\underline x} + \sum_{k=0}^{K_j-1} \alpha_{k,j} \one_{(x_{k,j}, x_{k+1,j})} +  \alpha_{K_j,j}\delta_{\overline x} .
\end{displaymath}
The $d=K_1+K_2+4$ coefficients $ \alpha_{k,j}$ are nonnegative and such that 
\begin{displaymath}
  \alpha_{-1,1}+  \alpha_{-1,2}+ \sum_{j=1}^2 \sum_{k=0}^{K_j-1} \alpha_{k,j} h_{k,j}+    \alpha_{K_1,1}+  \alpha_{K_2,2}=1.
\end{displaymath}
Let $\cP$ be the map the vector $\R^d\ni ( \alpha_{-1,1},\dots,   \alpha_{K_1,1}, \alpha_{-1,2},\dots,   \alpha_{K_2,2}) \mapsto m$. 
\subsubsection{Transport of $m$}
Let us set 
\begin{equation}
    \begin{aligned}[t]
  \chi_{-1,j}(x,y)&= \one\left(
      x=\underline x \hbox{  and }
     y = y_j
\right),
\\
  \chi_{k,j}(x,y)&= \one\left(
      x\in   (x_{k,j}, x_{k+1,j}] \hbox{  and }
     y = y_j
\right), \quad\quad \hbox{ if }  k=0,\dots,K_j-1, \\
  \chi_{K_j,j}(x,y)&= \one\left(
      x> \overline x  \hbox{  and }
     y = y_j
\right).
\end{aligned}
\end{equation}
The transported measure $m^*$ is approximated as follows: we fix a large integer $N$ and set $x_{n,k,j}= x_{k,j} +\frac n {N+1} h_{k,j}$ for $n\in \{1,\dots,N\}$,
$j=1,2$, and $k\in \{0,\dots, K_j-1\}$. Each point $x_{n,k,j}$ is transported by the optimal strategy to 
$\widehat x_{n,k,j}= x_{n,k,j}+ \Delta t s^*_{i,j} ( x_{n,k,j}, m)$. We then draw $\widehat y_{n,k,j}$:
\begin{displaymath}
  \PP(\widehat y_{n,k,j}=y_j)=1-\lambda_j \Delta t ,\quad \quad   \PP(\widehat y_{n,k,j}=y_{\bar \jmath})=\lambda_j \Delta t .
\end{displaymath}
Similarly,
\begin{itemize}
\item  $\underline x$ is transported by the optimal strategy to 
$\widehat x_{-1,j}= \underline x+ \Delta t s^*_{i,j} ( \underline x, m)$. We draw $\widehat y_{-1,j}= y_j $ or $=y_{\bar \jmath}$ with respective probabilities $1-\lambda_j \Delta t$ and $\lambda_j \Delta t$
\item  $\overline x$ is transported by the optimal strategy to 
$\widehat x_{K_j,j}= \overline x+ \Delta t s^*_{i,j} ( \overline x, m)$. We draw $\widehat y_{K_j,j}= y_j $ or $=y_{\bar \jmath}$ with respective probabilities $1-\lambda_j \Delta t$ and $\lambda_j \Delta t$.
\end{itemize}
The random variables  $\widehat y_{-1,j},\;\widehat y_{n,k,j}, \;\widehat y_{K_j,j}$ are all independent.

The new probability measure $m^* (\Delta t) $ is given by $m^*(\Delta t) =\sum_{j=1}^2 m^* _j \otimes \delta_{y=y_j} $, where
\begin{displaymath}
  m^*_j=  \alpha^* _{-1,j}\delta_{\underline x} + \sum_{k=0}^{K_j-1} \alpha^*_{k,j} \one_{(x_{k,j}, x_{k+1,j})} +  \alpha^*_{K_j,j}\delta_{\overline x},
\end{displaymath}
with
\begin{equation*}
  \alpha^*_{k,j}= \sum_{\ell=1} ^2 \left(  \begin{array}[c]{ll}
    & \ds \alpha_{-1,\ell} 
\chi_{k,j}\left( \widehat x_{-1,\ell}, \widehat y_{-1,\ell} \right)\\
&\ds + \frac 1 N  \sum_{p=0}^{K_\ell-1} \sum_{n=1}^{N}\alpha_{p,\ell}  \chi_{k,j}\left( \widehat x_{n,p,\ell} , \widehat y_{n,p,\ell}\right)\\
&\ds + \alpha_{K_\ell,\ell}  \chi_{k,j}\left( \widehat x_{K_\ell,\ell}, \widehat y_{K_\ell,\ell} \right)
                                           \end{array}
                                         \right).
\end{equation*}

\subsection{Approximation of the value functions $v_{i,j}$ with neural networks and a fixed point strategy}
\subsubsection{Approximation of the value functions}
Recall that the dimension of the space of discrete probability measures on $[\underline x, \overline x] \times \{y_1,y_2\}$ is $d$,
and that $\cP(M)$ denotes the measure associated to $M\in \R^d$.
It may be convenient to let the approximate value functions $v_{i,j}$
actually depend on less than $1+d$ parameters. Indeed, Krusell and Smith have conjectured that
the value function $v_{i,j}$ mostly depends  on $m$ though the interest rate  given by \eqref{eq:5}.
We therefore introduce an integer $0\le d_0 < d $, and a map $\cF_{i,j}: M\in\R^d \mapsto \cF_{i,j}(M)\in \R^{d_0+1}$,  $\cF_{i,j}(M)$ being a collection of relevant parameters that can be constructed from the probability measure $\cP(M)$.
For example, such  parameters may include the interest rate $r$ given by \eqref{eq:5} and some moments of $m_{k}$, $k=1,2$.
\\
We are going to approximate the value functions $v_{i,j}$ by means of neural networks,
exploiting their capability to provide appropriate sets of parameterized functions.
In our strategy, the first component of $\cF_{i,j}$ is the interest rate, while
the last $d_0$ components of $\cF_{i,j}$ are not determined beforehand, but are rather found in an adaptive manner as an output of the first layer of a neural network. 
\\
Let $\cN_{L,d,d_0, d_1,\dots, d_L}$ denote the chosen set of neural networks, which are real valued functions  with the following characteristics: the number of layers is $L+1$, the input dimension is $d+1$ where $d$ has been introduced above, the output dimension $d_L$ is $1$,
and the number of neurons in the hidden layers are $d_1,\dots, d_{L-1}$. We  approximate $  v_{i,j}(x,\cP(M)) $ by $ v^{\cN}_{ i,j}  (x,\cP(M))=   \cV_{i,j}  (x, \cF_{i,j}(M))$, where
$  \cV_{i,j}\in  \cN_{L,d,d_0, d_1,\dots, d_L}$.

\subsubsection{A fixed point strategy}
Before describing the fixed point strategy, let us introduce a large set $\cS$  of samples $  (x, M)\in  [\underline x, \overline x] \times \R^d$ such that $\cP(M)$ is a probability measure on $ [\underline x, \overline x] \times\{y_1,y_2\}$.
\\
We consider the following fixed point iterations: $(v^{\cN, n}_{ i,j})_{i,j}\to (v^{\cN, n+1}_{ i,j})_{i,j}$:

\begin{itemize}
\item For $i=1,2$, $j=1,2$, set $v_{i,j}=v^{\cN, n}_{ i,j}$
\item For $i=1,2$, $j=1,2$,
  \begin{itemize}
  \item 
 For each $(x,M)\in \cS$, set $m=\cP(M)$ and compute $c_{i,j}^* (x,m)$ by \eqref{eq:7}.
\item Find $v^{\cN, n+1}_{ i,j}$ as the minimizer in the class of functions described above, of
  \begin{displaymath}
   V \mapsto  \sum_{(x,M)\in \cS} \left| R_{i,j}(x,m)\right|^2
  \end{displaymath}
  where

\begin{eqnarray*}
      R_{i,j}(x,m)
      = &&
           (1+\rho \Delta t ) V (x,m)
\\
          && - v_{i,j}\Bigl(x+ \Delta t (w_i(m)y_j + r_i(m) x -c^*_{i,j}(x,m))  , m^*(\Delta t)    \Bigr)\\
&&+\lambda_j \Delta t \left(  v_{i,j}(x,m)  -  v_{i,\bar \jmath}(x,m) \right)\\ && + \mu_i \Delta t \left(  v_{i,j}(x,m)  -  v_{\bar  \imath,j}(x,m) \right) +\Delta t u \left(c^*_{i,j}(x,m)\right)
  \end{eqnarray*}
  and  $m^* (\Delta t) $ is the transported version of $m$ computed as in Subsection \ref{sec:approximation-m-its}.
\end{itemize}
  \end{itemize}

  \section{Some results}

  Hereafter, we present prelimininary numerical results. We insist that these results are only  meant to illustrate the method and its outputs.
  We need to run more simulations, on a larger scale,  to be more confident on the results and draw sound conclusions.

We took $\underline x= 0$, $\overline x= 30$, $y_1=0.7$, $y_2=1.4$, $\lambda_1=0.05$, $\lambda_2=0.1$, $\rho=0.15$, $\delta=0.05$, $\alpha=0.5$, $A_1=0.9$, $A_2=1.1$, $\mu_1=\mu_2=0.2$.
\\
We used a time step of $0.25$ year.
\\
The grids used to discretize the densities $m_0$ and $m_1$ have  $d=17+10=27$ nodes. The repartition of the nodes
is chosen in such a way that the grid steps are equally weighted by the  measure found at the equilibrium in Aiyagari's model with $A=1$.

.

\subsection{An  architecture designed  for exploration}

To construct  approximations of thes solutions, we need to design a neural network architecture providing parameterized functions and remaining consistent with the economics of Krusell-Smith model.

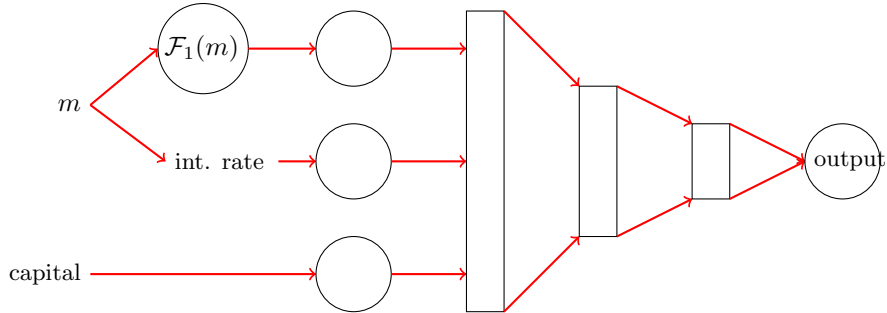
\begin{figure}[H]
  \begin{center}
    \begin{tikzpicture}[scale=0.5, trans/.style={thick,<->,shorten >=2pt,shorten <=2pt,>=stealth} ]
      \draw(0,4.5) node[left]{{ $m$}};
      \draw (3,6) circle (1.2) ;
       \draw(4.2,6) node[left]{{$\cF_1(m)$}};
       \draw[->,red,thick] (0,4.5) -- (2,3);
        \draw(2,3) node[right]{{\small int. rate}};
      \draw[->,red,thick] (0,4.5) -- (1.8,6);
      \draw(0,0) node[left]{{\small capital}};
      \draw[->,red,thick] (0,0) -- (6,0);
      \draw[->,red,thick] (4.2,6) -- (6,6);
      \draw (7,6) circle (1) ;
      \draw(6.5,6) node[right]{{}};
      \draw[->,red,thick] (5,3) -- (6,3);
      \draw (7,3) circle (1) ;
      \draw(6.5,3) node[right] {{}};
       \draw (7,0) circle (1) ;
       \draw(6.5,0) node[right] {{}};
       \draw[->,red,thick] (8,3) -- (10,3);
       \draw[->,red,thick] (8,6) -- (10,6);
        \draw[->,red,thick] (8,0) -- (10,0);
        \draw (10,-1) rectangle (11,7);
        \draw[->,red,thick] (11,-1) -- (13,1);
         \draw[->,red,thick] (11,7) -- (13,5);
         \draw (13,1) rectangle (14,5);
          \draw[->,red,thick] (14,1) -- (16,2);
         \draw[->,red,thick] (14,5) -- (16,4);
         \draw (16,2) rectangle (17,4);
           \draw[->,red,thick] (17,2) -- (19,3);
         \draw[->,red,thick] (17,4) -- (19,3);
         \draw (20,3) circle (1) ;
           \draw(19,3) node[right] {{\small output}};
         \end{tikzpicture}
         \caption{The chosen architecture with $d_0=1$: note that in the first layer, the input $m$ is first mapped to the interest rate and to another variable $\cF_1(m)$. Krusell and Smith conjecture that the optimal value  depends on the distribution of states mostly through  the interest rate, so we expect that the auxiliary variable  $\cF_1(m)$ should have a rather small importance. In the simulations reported below, the architecture is as follows: all the layers or sublayers described below, except the final one, involve the {\sl softplus} activation function. The vector $(m, x)$  is processed by the first layer as follows: first $(m,x)$ is mapped to $(\cF_1(m), r, x)$ and the map $\cF_{1}$ is described by a first sublayer. This means that $\cF_1$ is the composition of an affine map from $\R^{d}$ to  $\R$  with the {\sl softplus} activation function.  The variable $\cF_1(m)$ is  mapped to a vector in  $\R^{80}$ by a second sublayer. The interest rate $r$ is  mapped to a vector in  $\R^{20}$ by a third sublayer.
The capital $x$  is mapped to a vector in $\R ^ {150}$ by a fourth sublayer. The output of the first layer is therefore  a vector in  $\R^{d_1}$,   $d_1=80+20+150=250$, obtained by concatenation. This information is then processed by a sequence of $4$ hidden fully connected layers, with dimensions $300$, $150$, $50$ and $20$.  }
 \label{fig:0}
  \end{center}
\end{figure}

Our goal is to find neural networks approximations of the value function $v_{i,j}(x,m)$  in the four different situations indexed by $i=1,2$, $j=1,2$ corresponding to  productive/unproductive households, slow/fast economy i.e. $A=A_1,A_2$,
as a function of $x, r, \cF_{1,i,j}(m), \cdots,\cF_{d_0,i,j}(m)$, where $\cF_{1,i,j}(m), \cdots,$ $\cF_{d_0,i,j}(m)$ are found in adaptative way as the outputs of sublayers contained in the first layer in the neural network.
In our strategy, we aim at starting with $d_0=0$, then increasing $d_0$ gradually. For example, with $d_0=1$, the neural network architecture is meant to find $\cF_{1,i,j}(m)$ in order to complement the information given by the interest rate $r$.
In the simulations reported below, the four neural networks have all the same architecture,  displayed on Figure \ref{fig:0}.
Here, for brevity, we  will sometimes omit the indices $i$ and $j$, i.e. we will use the notation 
$\cF_1(m)$ for the  output of the first sublayer contained in  the first layer, see Figure \ref{fig:0}, remembering
that it will vary according to the considered situation, (productive/unproductive households , fast/slow economy).


\subsubsection{The optimal savings as a function of the capital and the interest rate}

In Figure \ref{fig:1}, we fix the value $\cF_{1,i,j}(m)$, and plot the contours of the optimal savings policy as a a function of the capital and the interest rate. The dotted lines correspond to negative values.

\begin{center}
  \includegraphics[width=0.4\linewidth]{{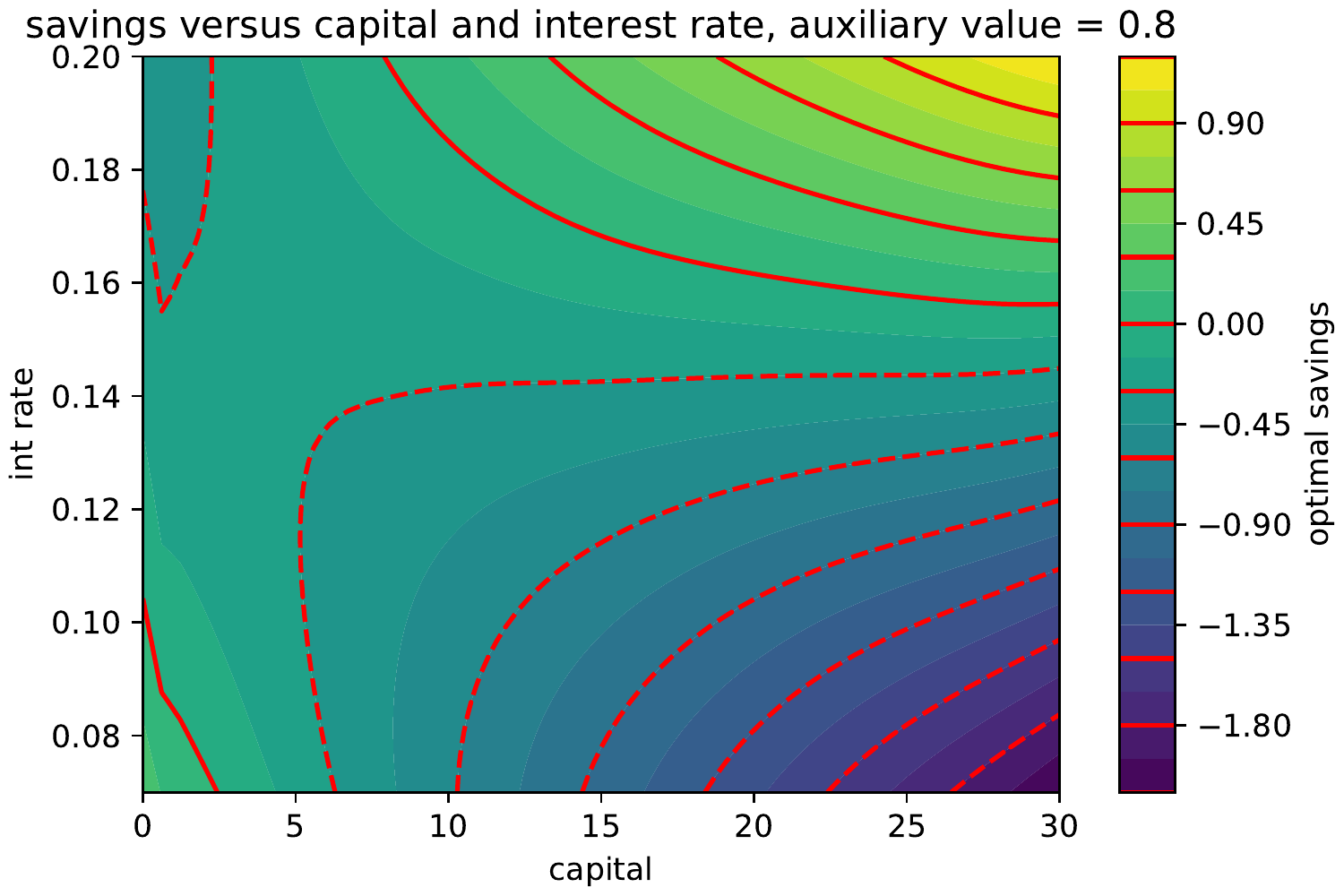}}
  \includegraphics[width=0.4\linewidth]{{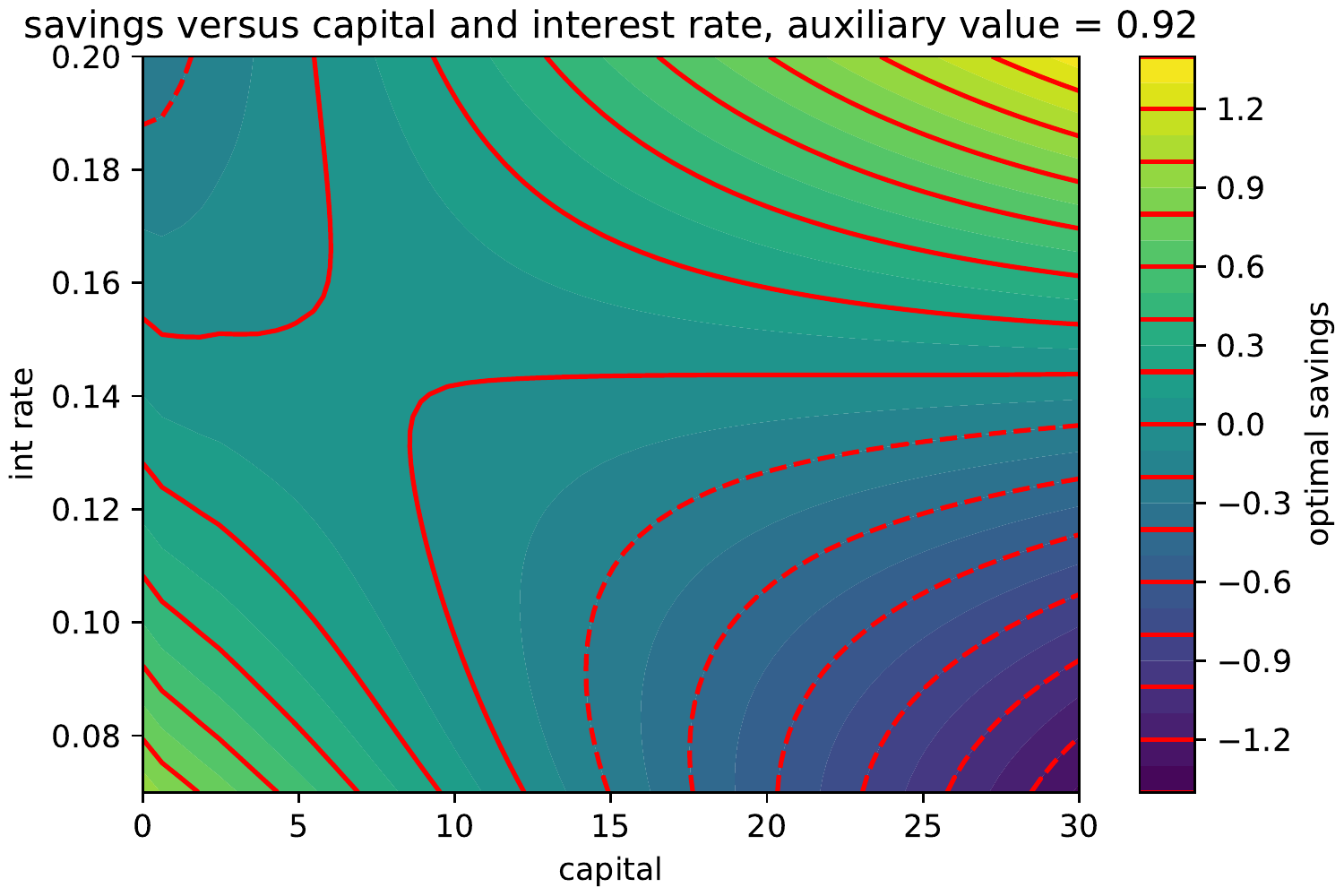}}
   \includegraphics[width=0.4\linewidth]{{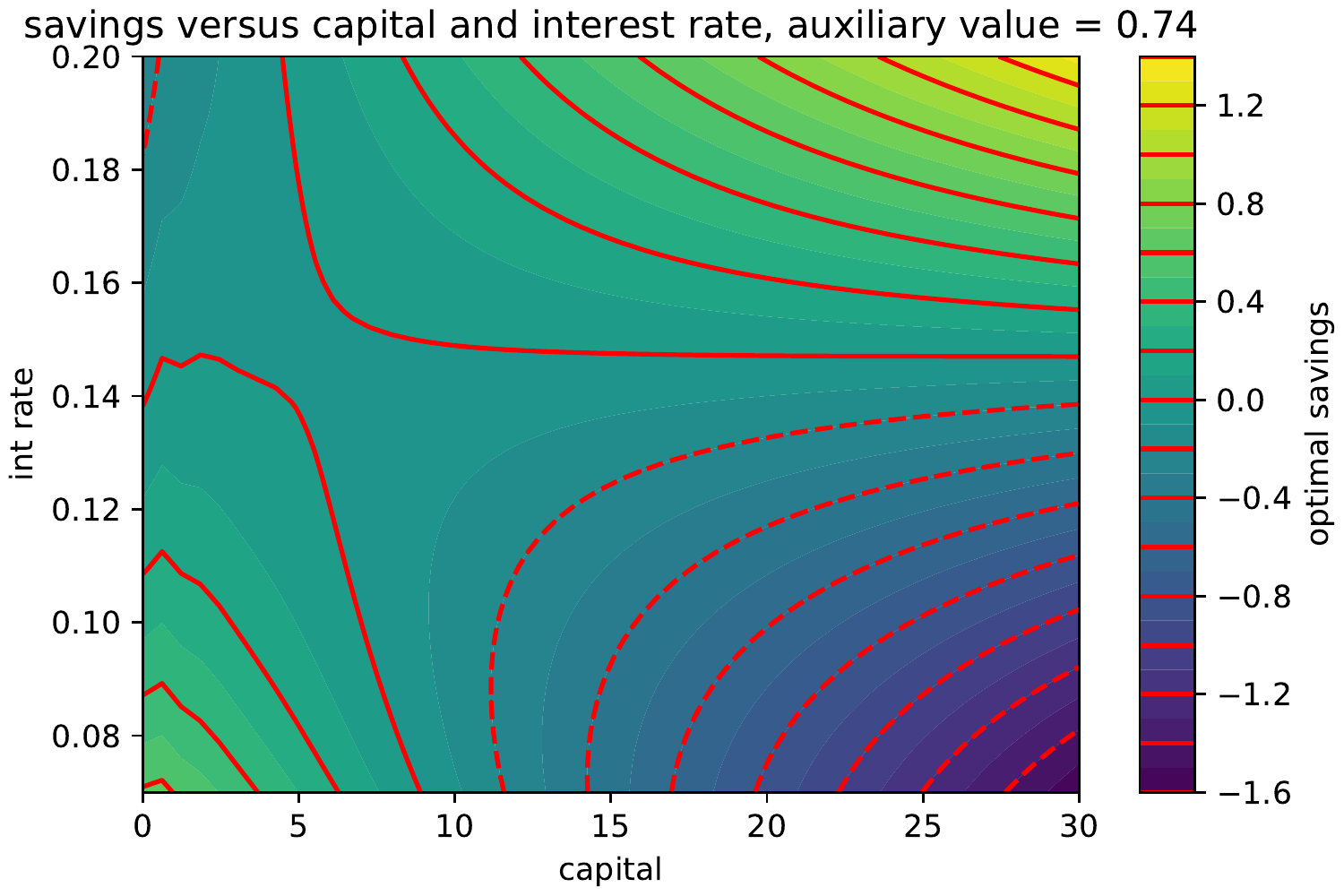}}
  \includegraphics[width=0.4\linewidth]{{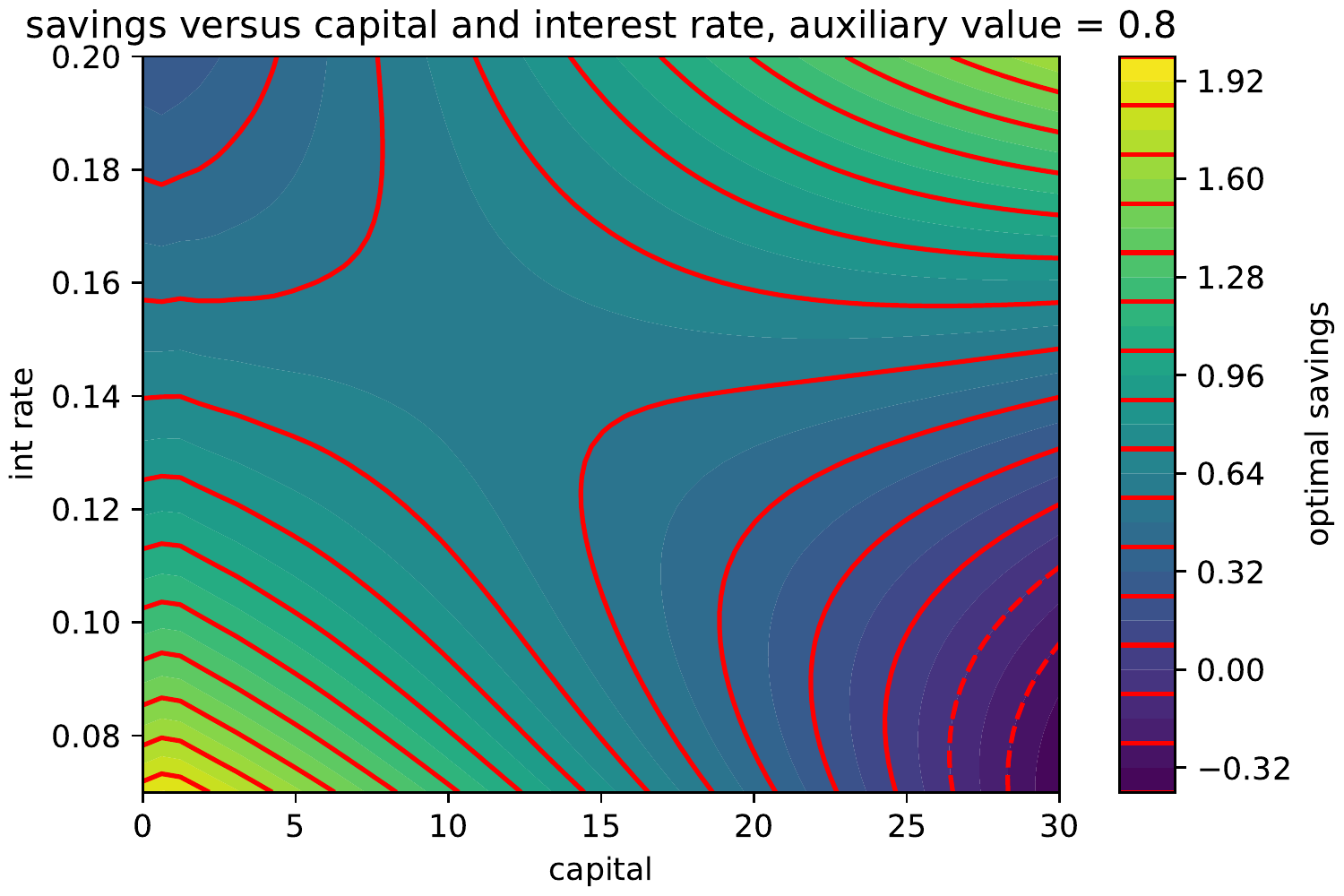}}
 \captionof{figure}{Top: Optimal savings for nonproductive(left) and productive(right) households when $A=A_1$  (slow economy) as a function of $x$ and $r$. Bottom:  Optimal savings for nonproductive(left) and productive(right) households when $A=A_2$ (fast economy).
 }
  \label{fig:1}
\end{center}

\subsubsection{The optimal savings as a function of the capital for a  given distribution $m$}
\begin{center}
  \includegraphics[width=0.4\linewidth]{{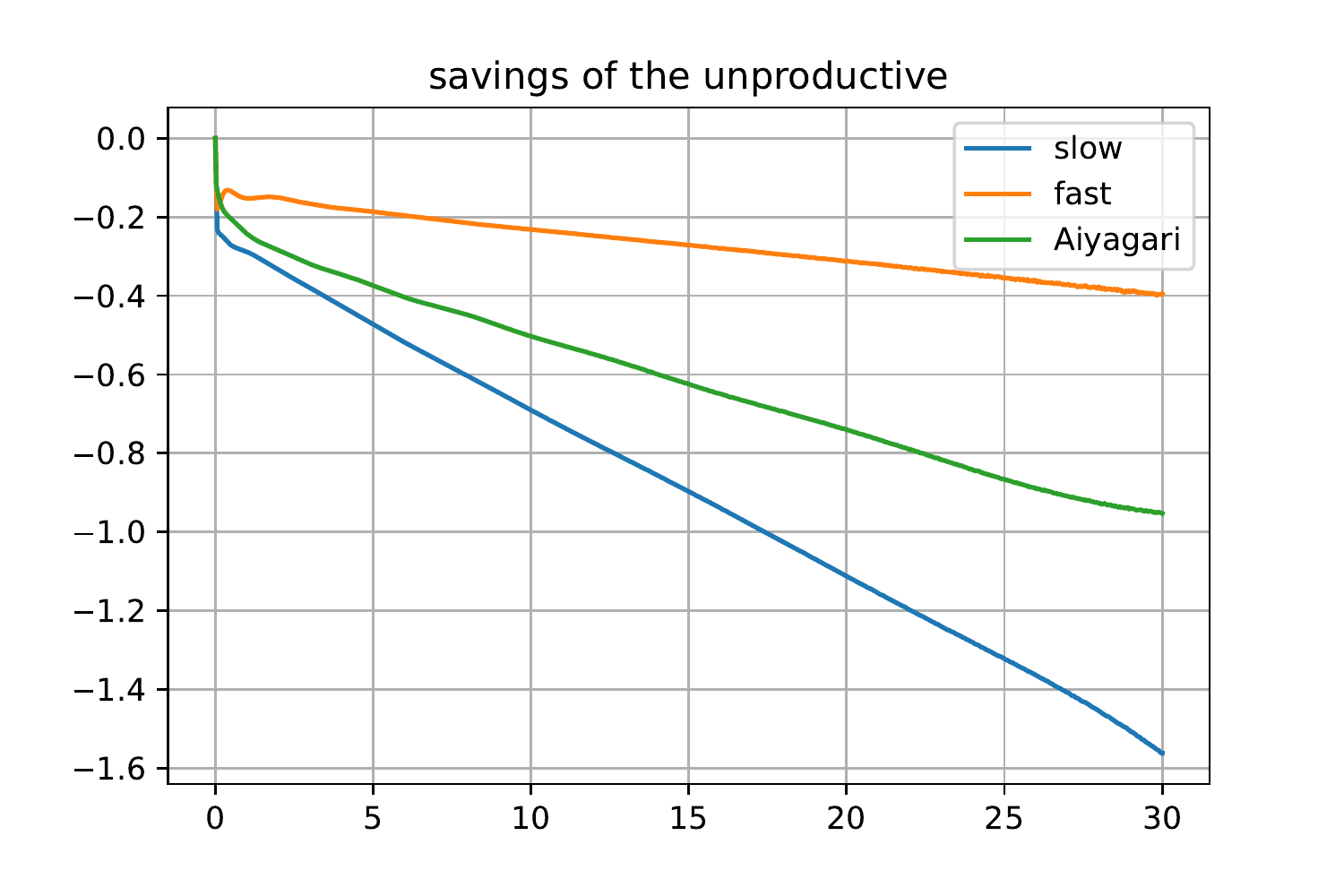}}
  \includegraphics[width=0.4\linewidth]{{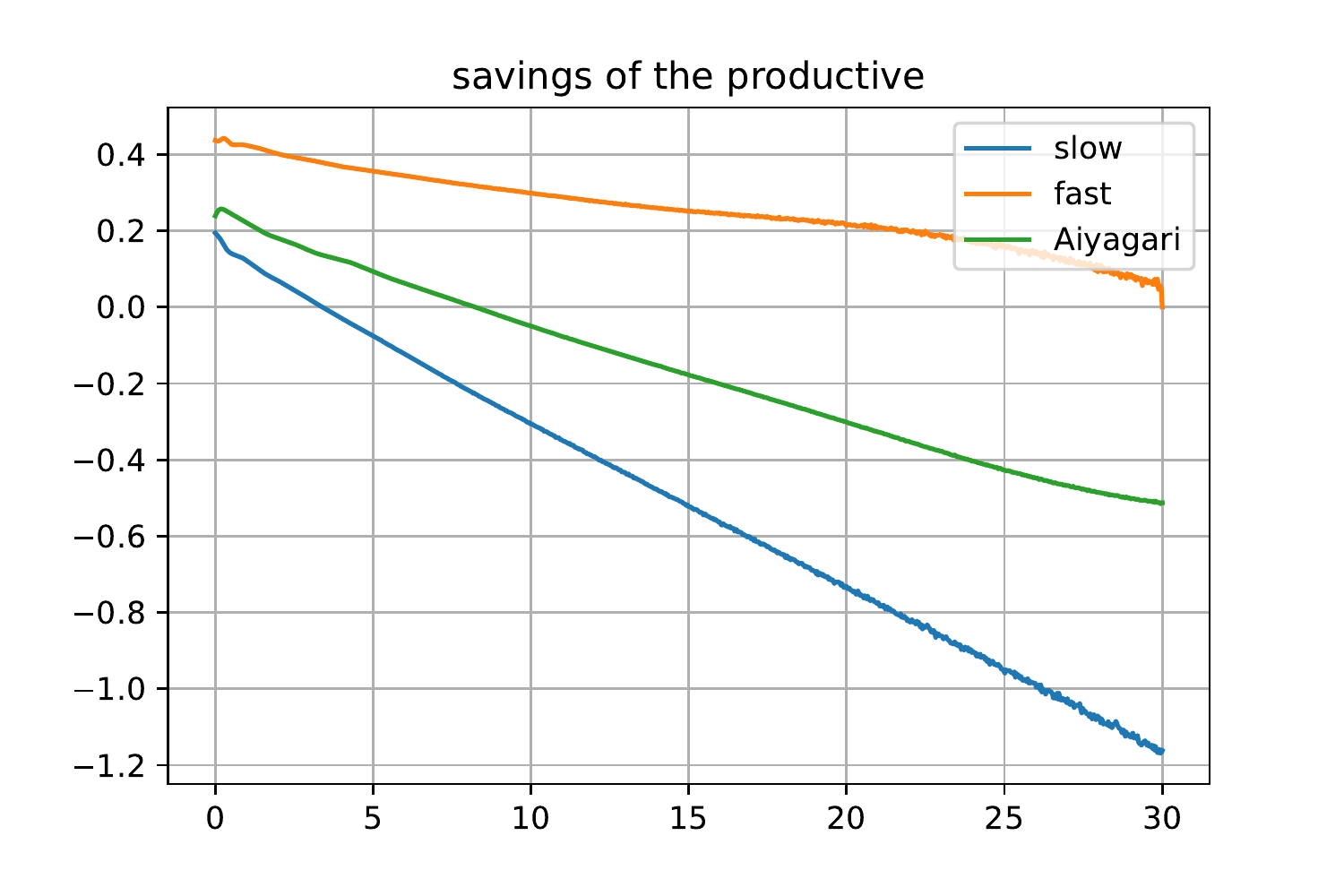}}
  \captionof{figure}{Optimal savings for nonproductive(left) and productive(right) households as a function of $x$ for the two values of $A$ (slow and fast economy), and the savings found by Aiyagari's model with $A=1$.
    The measure $m$ is chosen  as  the equilibrium in Aiyagari's model with $A=1$. Note the singular behavior in the savings of the unproductive near $x=0$; indeed, the credit constraint is biding in this situation. Note also that the savings policy found in the Aiyagari's model with $A=1$ lies below (resp. above) the policy related to the fast (resp. slow) economy.
 }
  \label{fig:1b}
\end{center}

\subsubsection{The optimal savings as a function of the interest rate and  the auxiliary variable $\cF_1(m)$}

\begin{center}
  \includegraphics[width=0.4\linewidth]{{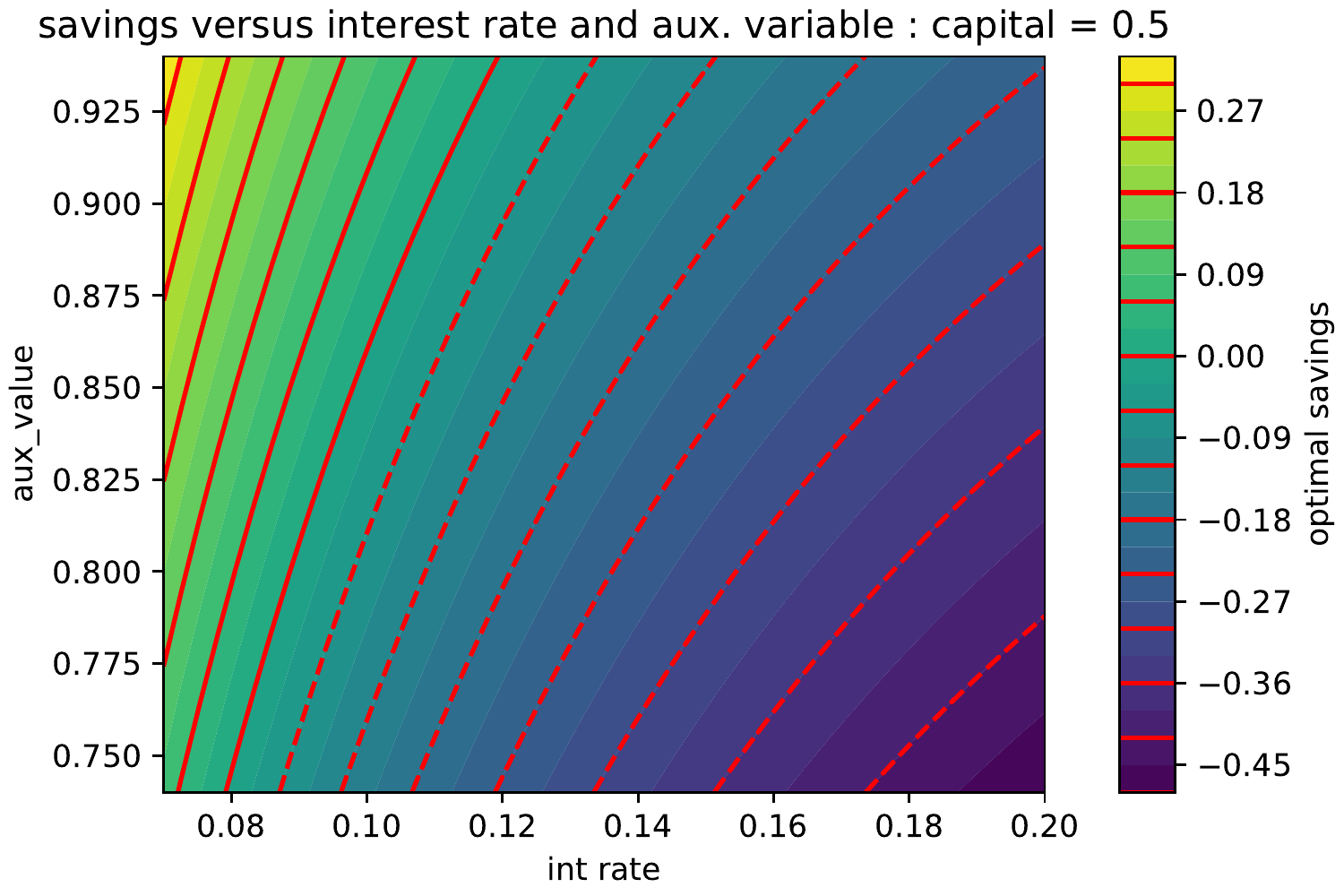}}
   \includegraphics[width=0.4\linewidth]{{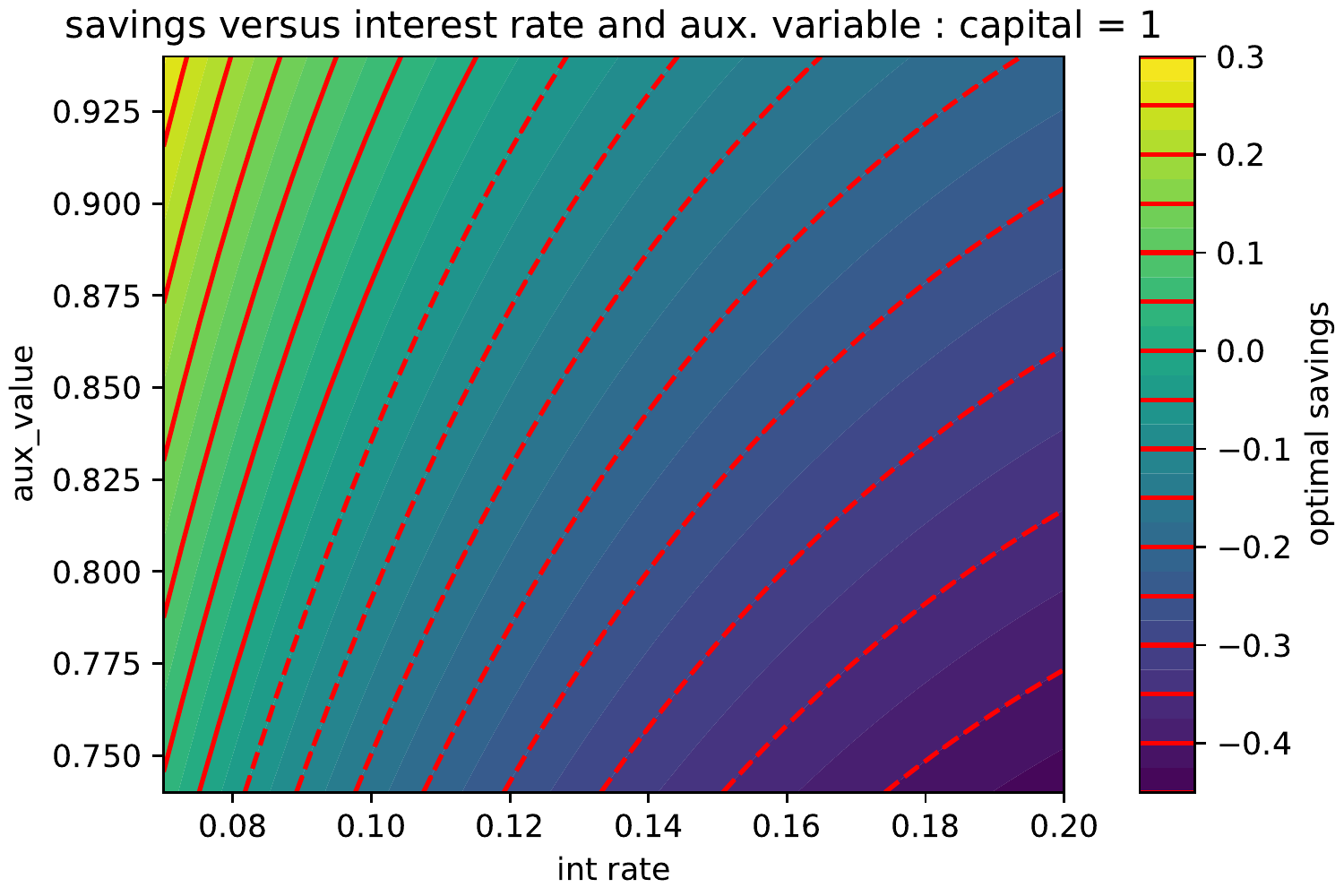}}
 \includegraphics[width=0.4\linewidth]{{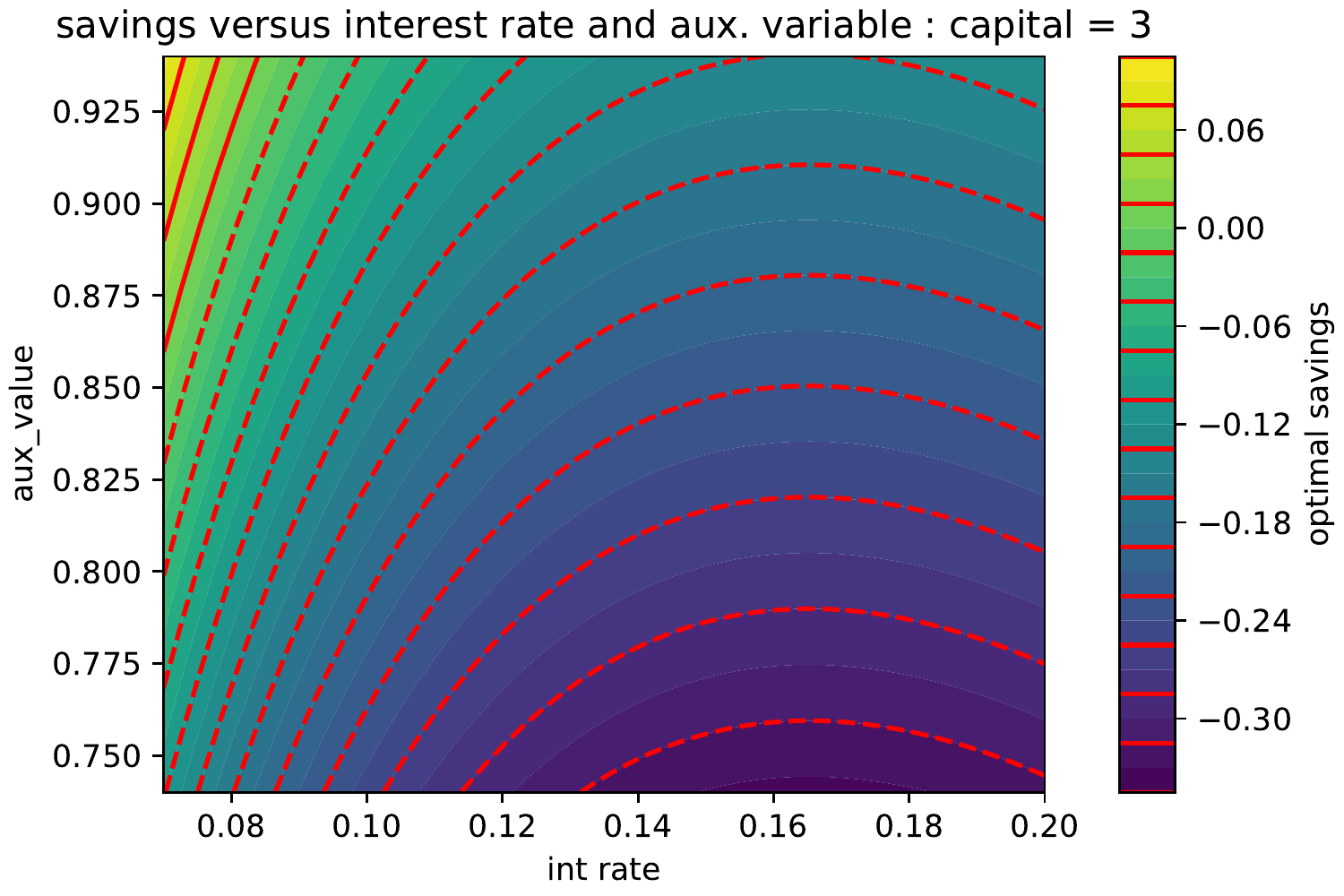}}
 \includegraphics[width=0.4\linewidth]{{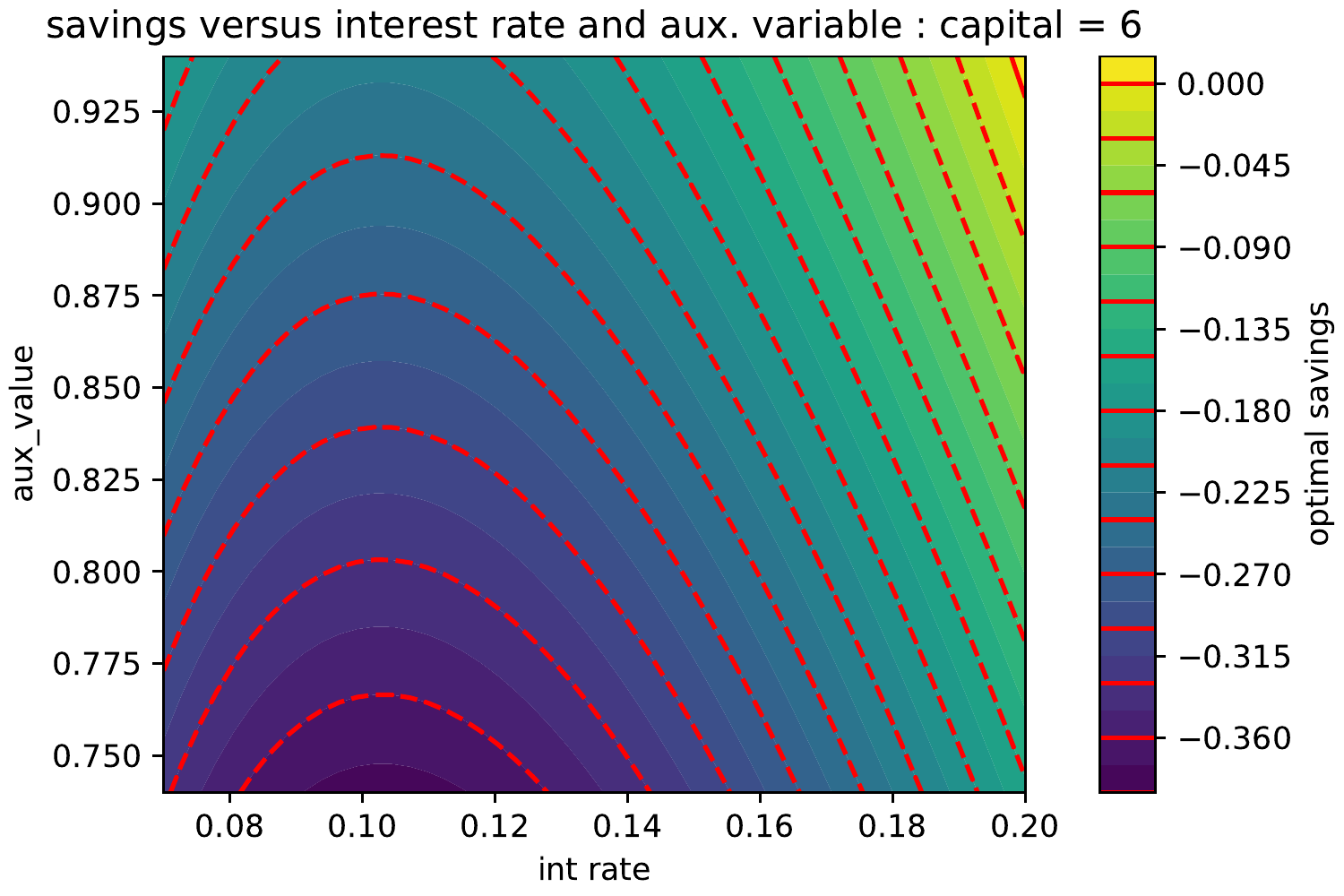}}
 \captionof{figure}{Optimal savings for nonproductive households when $A=A_1$  as a function of  $r$ and the auxiliary variable $\cF_{1}(m)$  for $x=0.5,\;1,\;3,\;6$.
 }
  \label{fig:2}
\end{center}

\begin{center}
 \includegraphics[width=0.4\linewidth]{{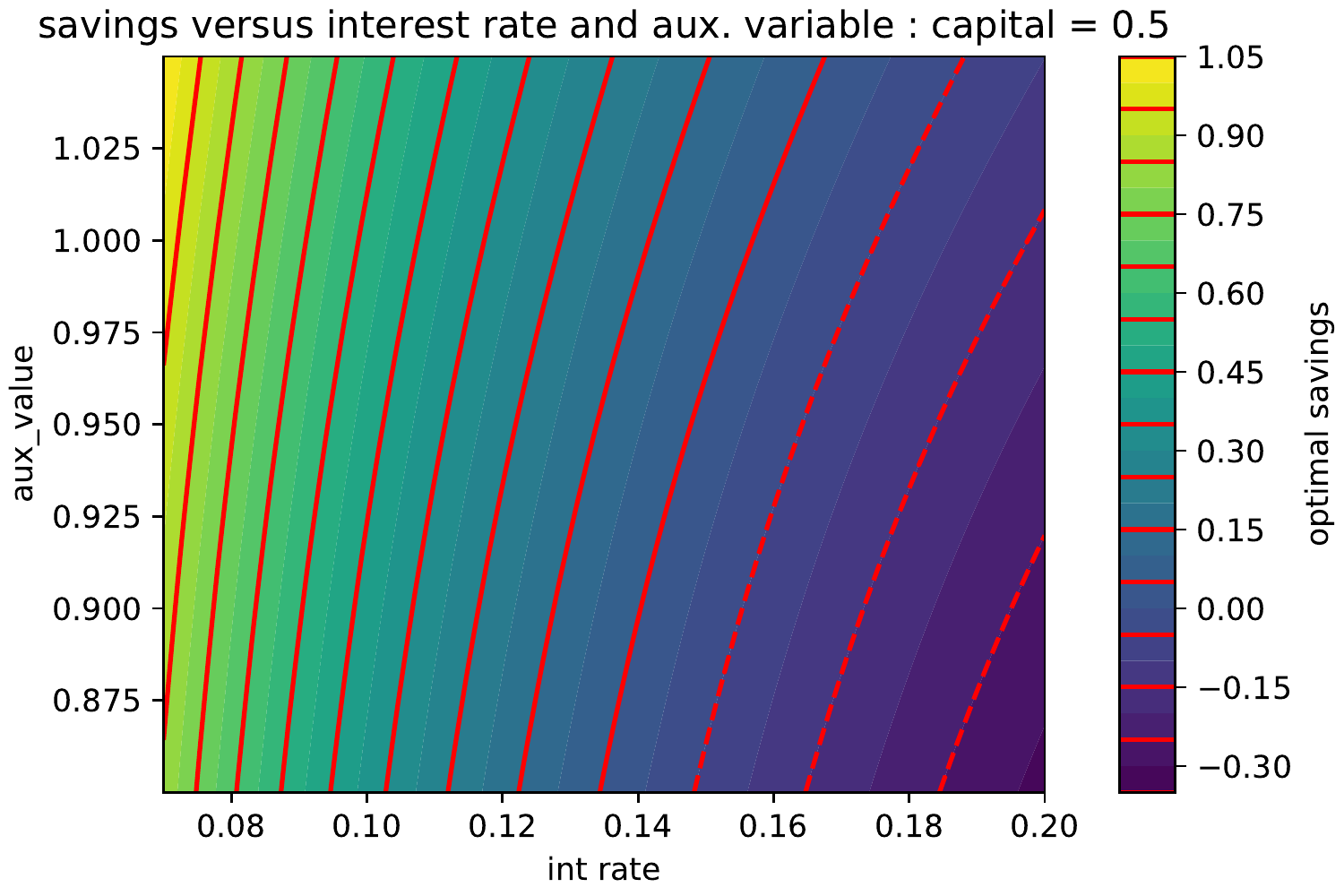}}
   \includegraphics[width=0.4\linewidth]{{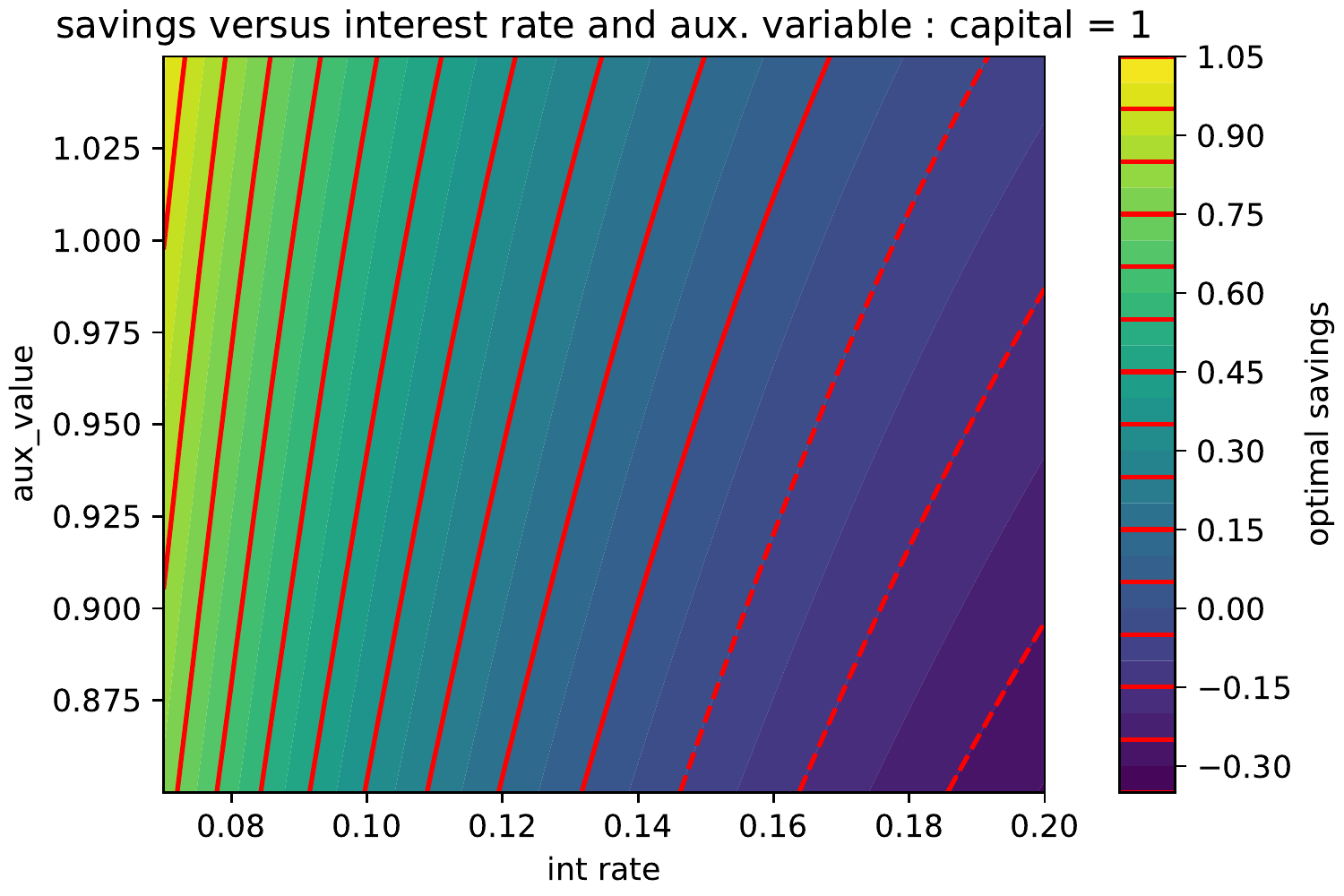}}
 \includegraphics[width=0.4\linewidth]{{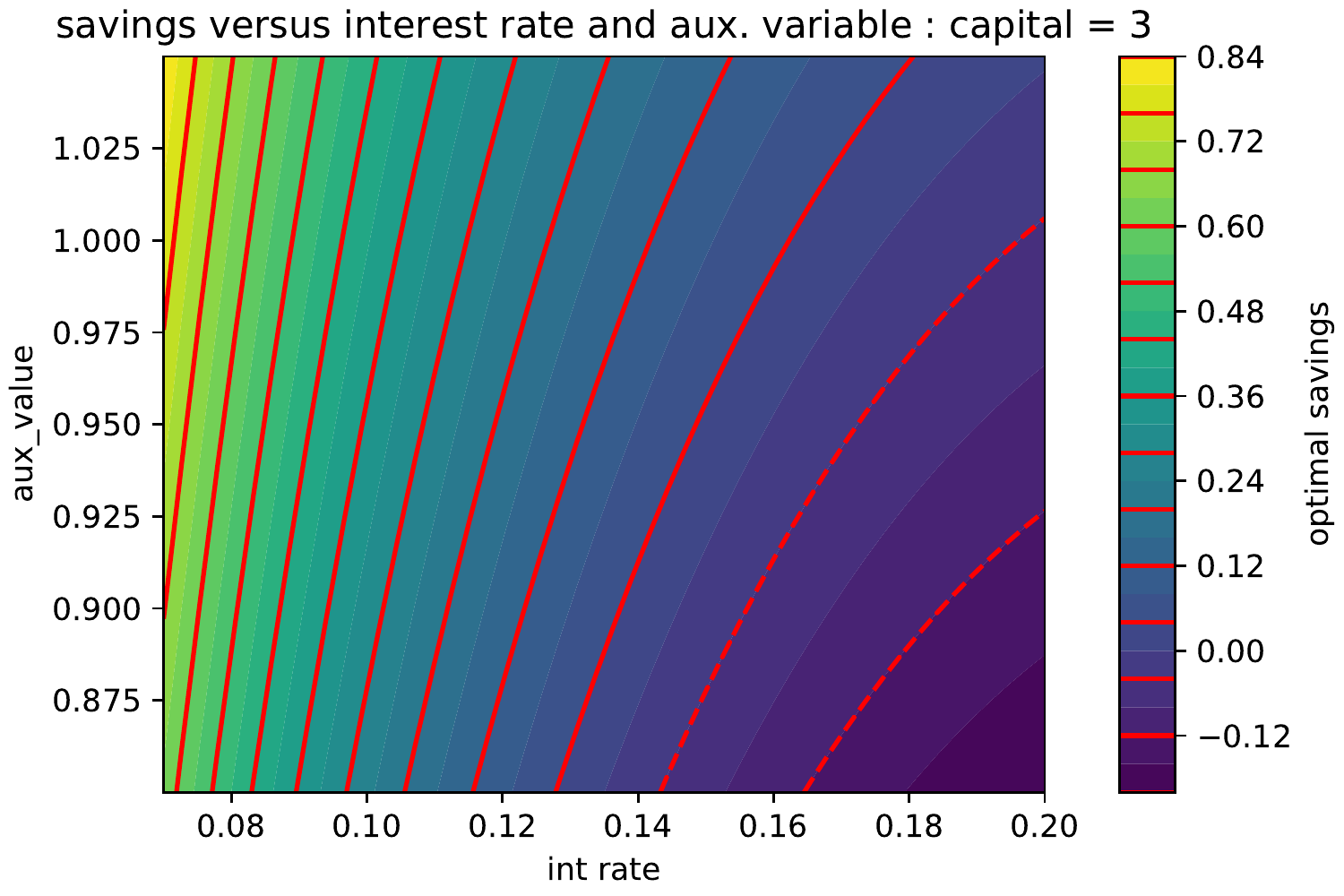}}
 \includegraphics[width=0.4\linewidth]{{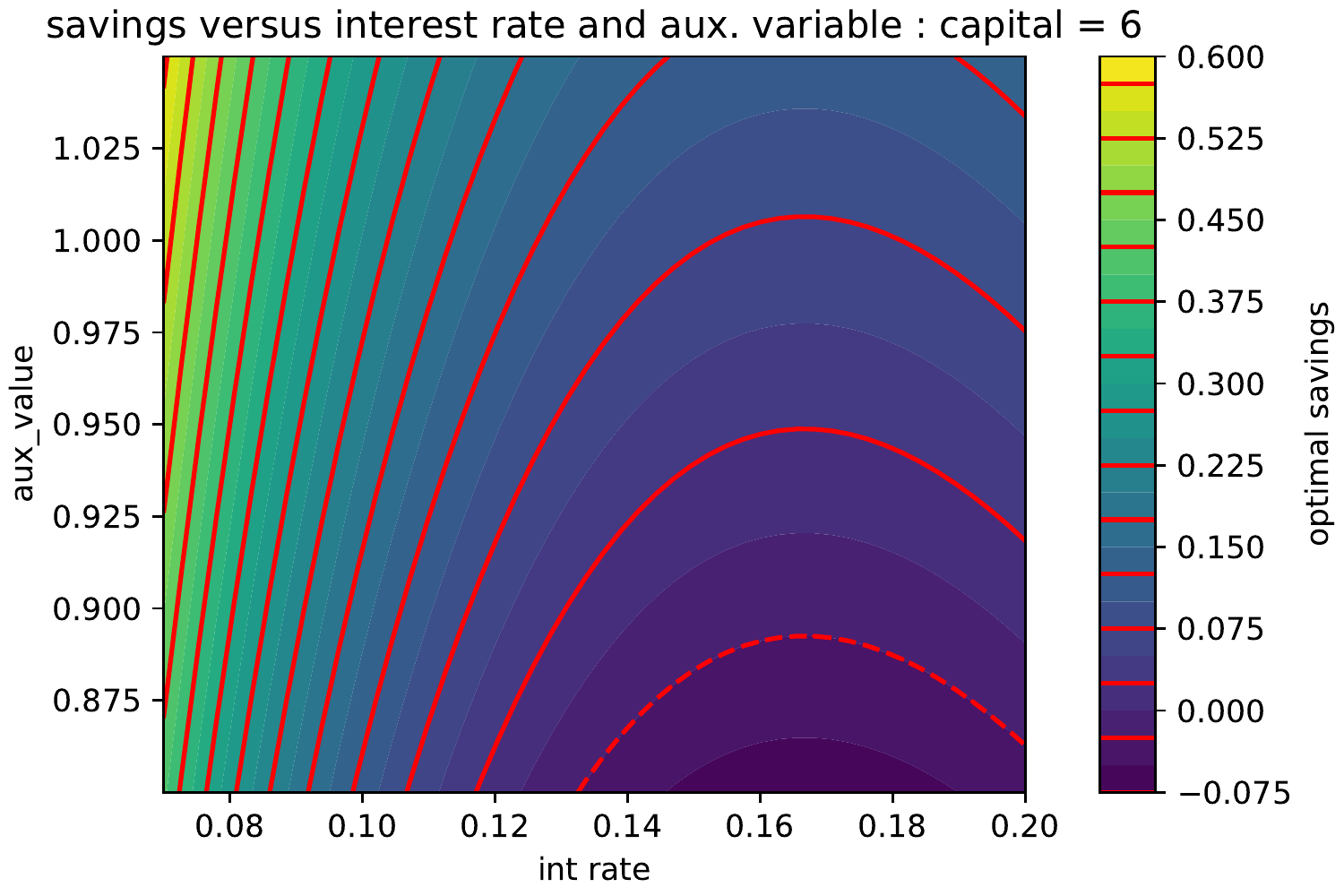}}
 \captionof{figure}{Optimal savings for productive households when $A=A_1$  as a function of  $r$ and the auxiliary variable $\cF_1(m)$ for $x=0.5,\;1,\;3,\;6$.
   }
  \label{fig:3}
\end{center}

\begin{center}
   \includegraphics[width=0.4\linewidth]{{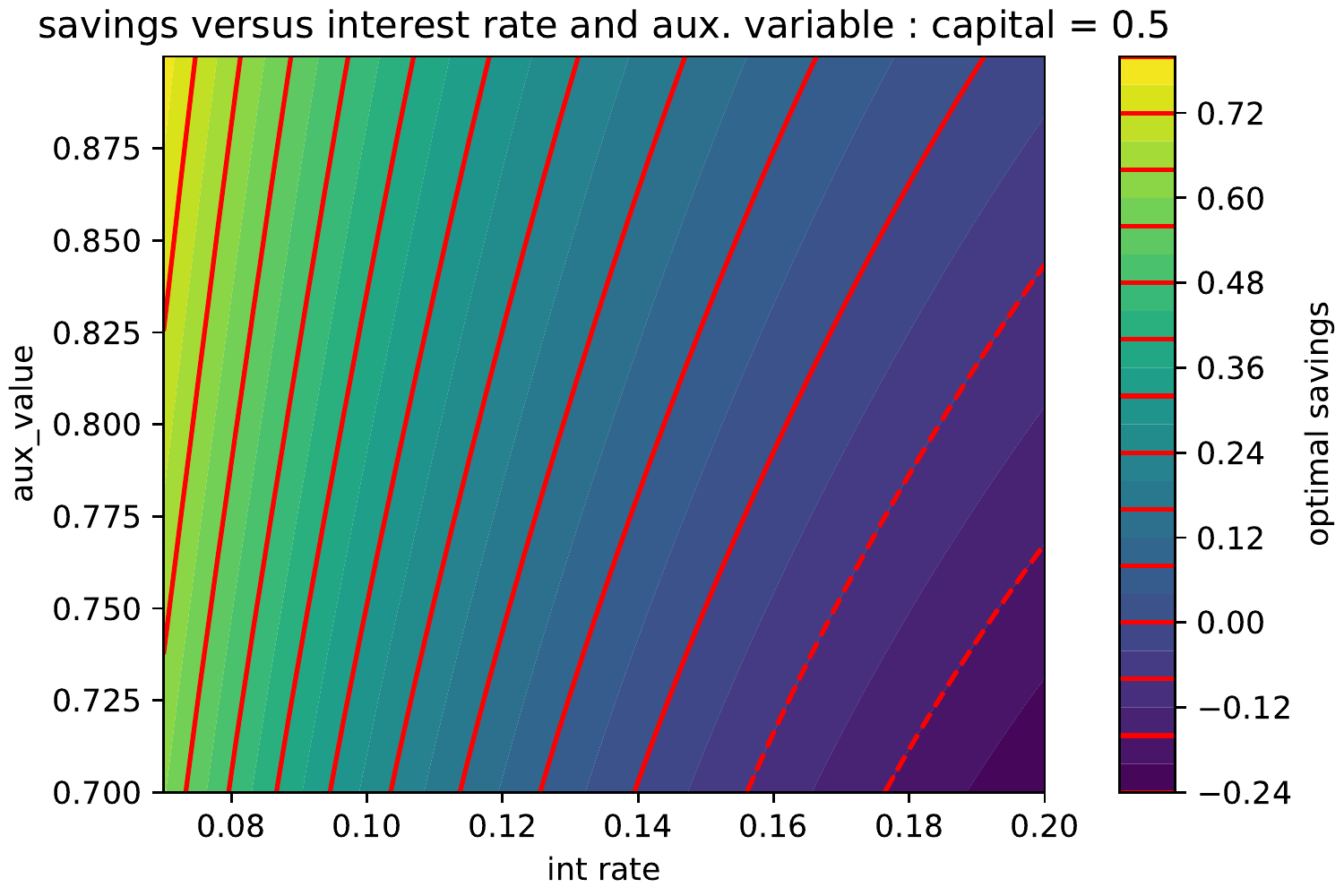}}
   \includegraphics[width=0.4\linewidth]{{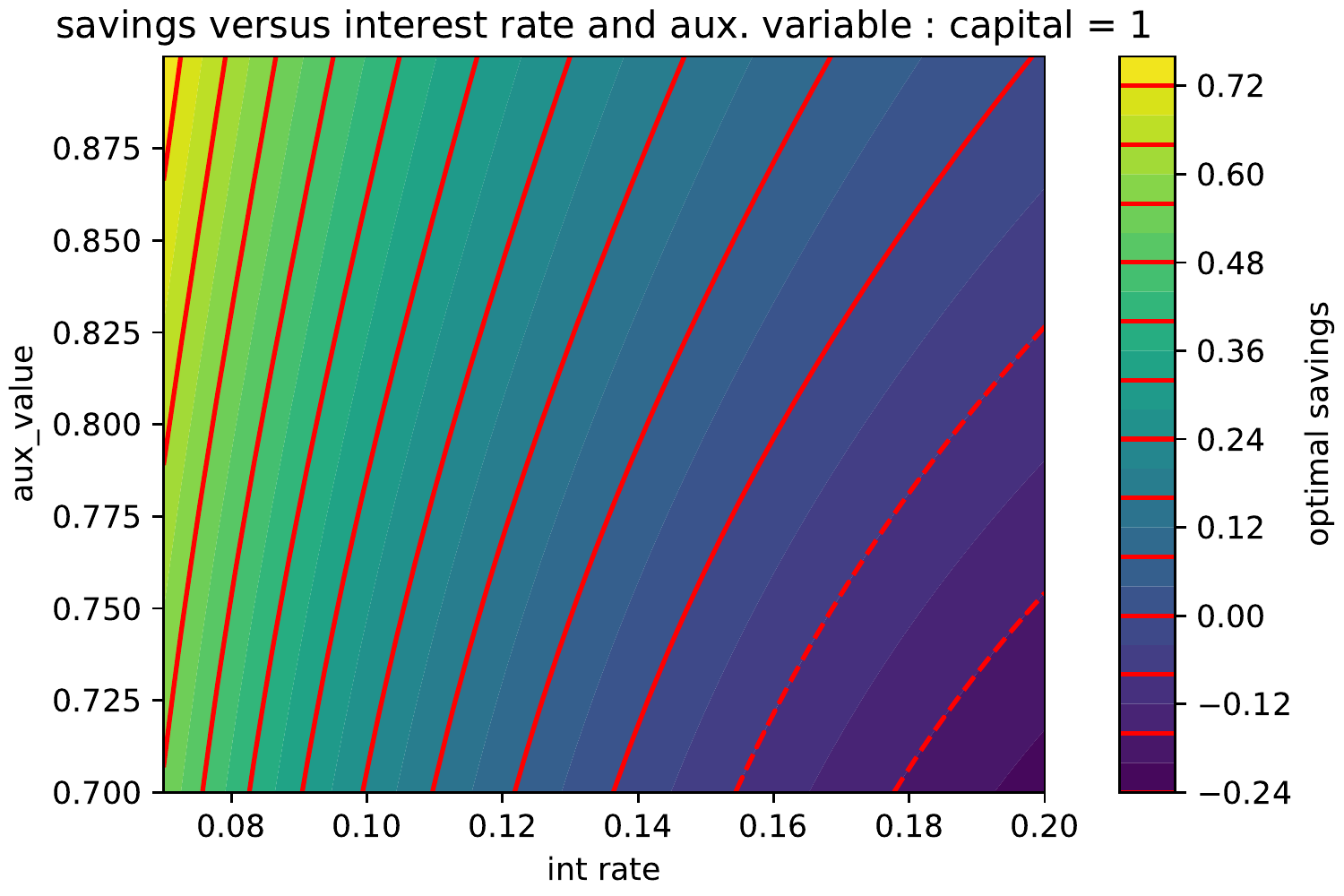}}
 \includegraphics[width=0.4\linewidth]{{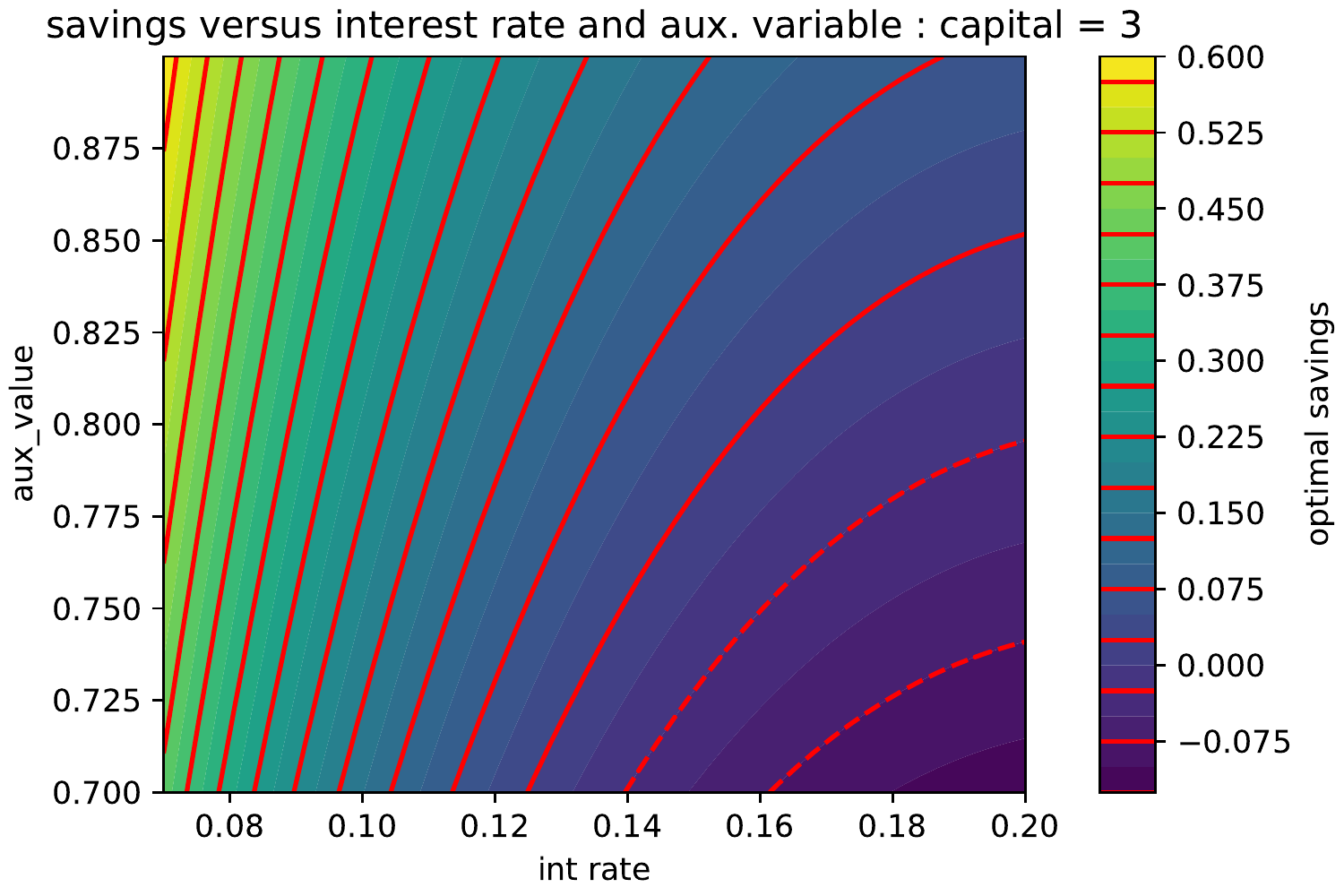}}
 \includegraphics[width=0.4\linewidth]{{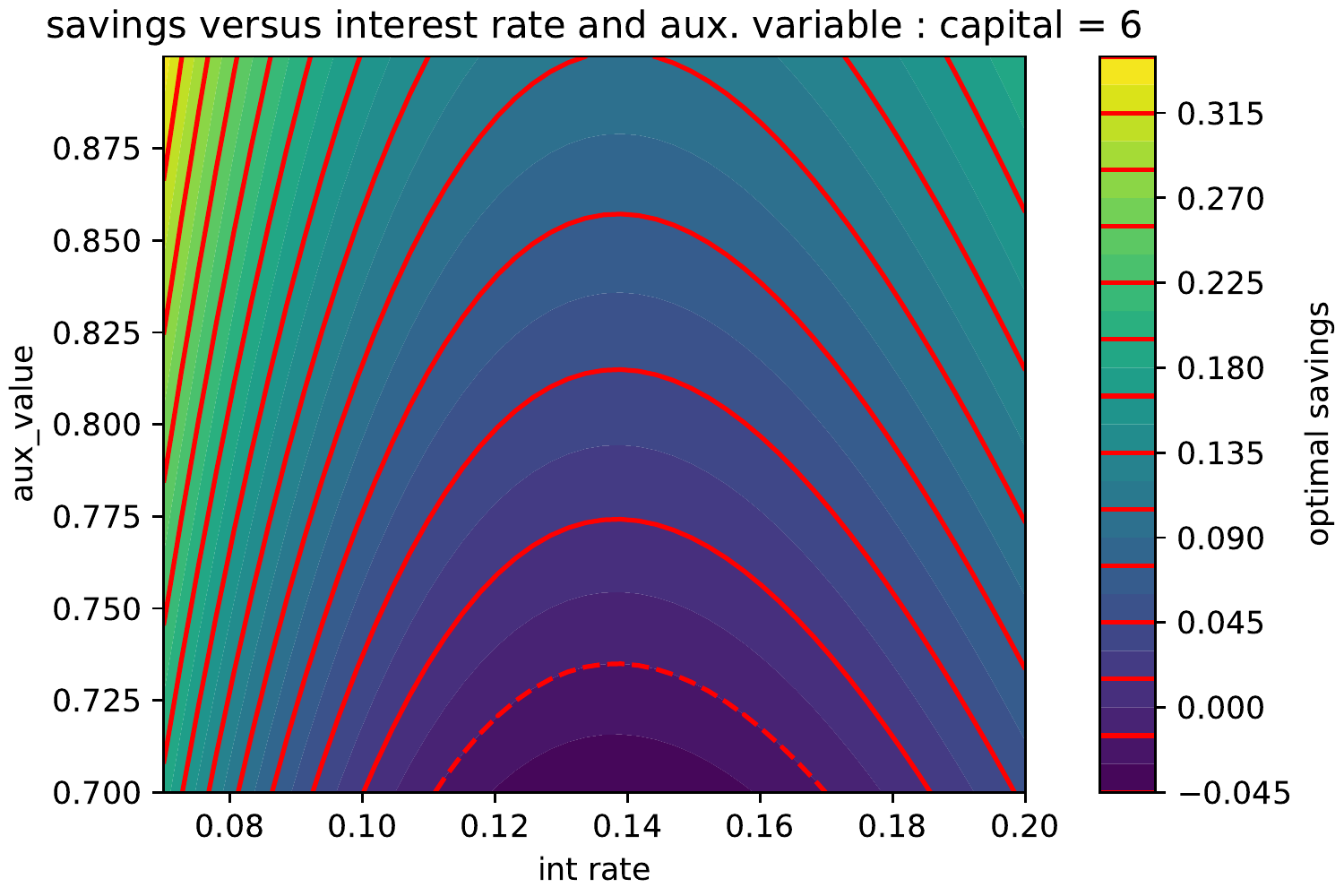}}
 \captionof{figure}{Optimal savings for nonproductive households when $A=A_2$  as a function of  $r$ and the auxiliary variable $\cF_1(m)$ for $x=0.5,\;1,\;3,\;6$. 
   }
  \label{fig:4}
\end{center}

\begin{center}
    \includegraphics[width=0.4\linewidth]{{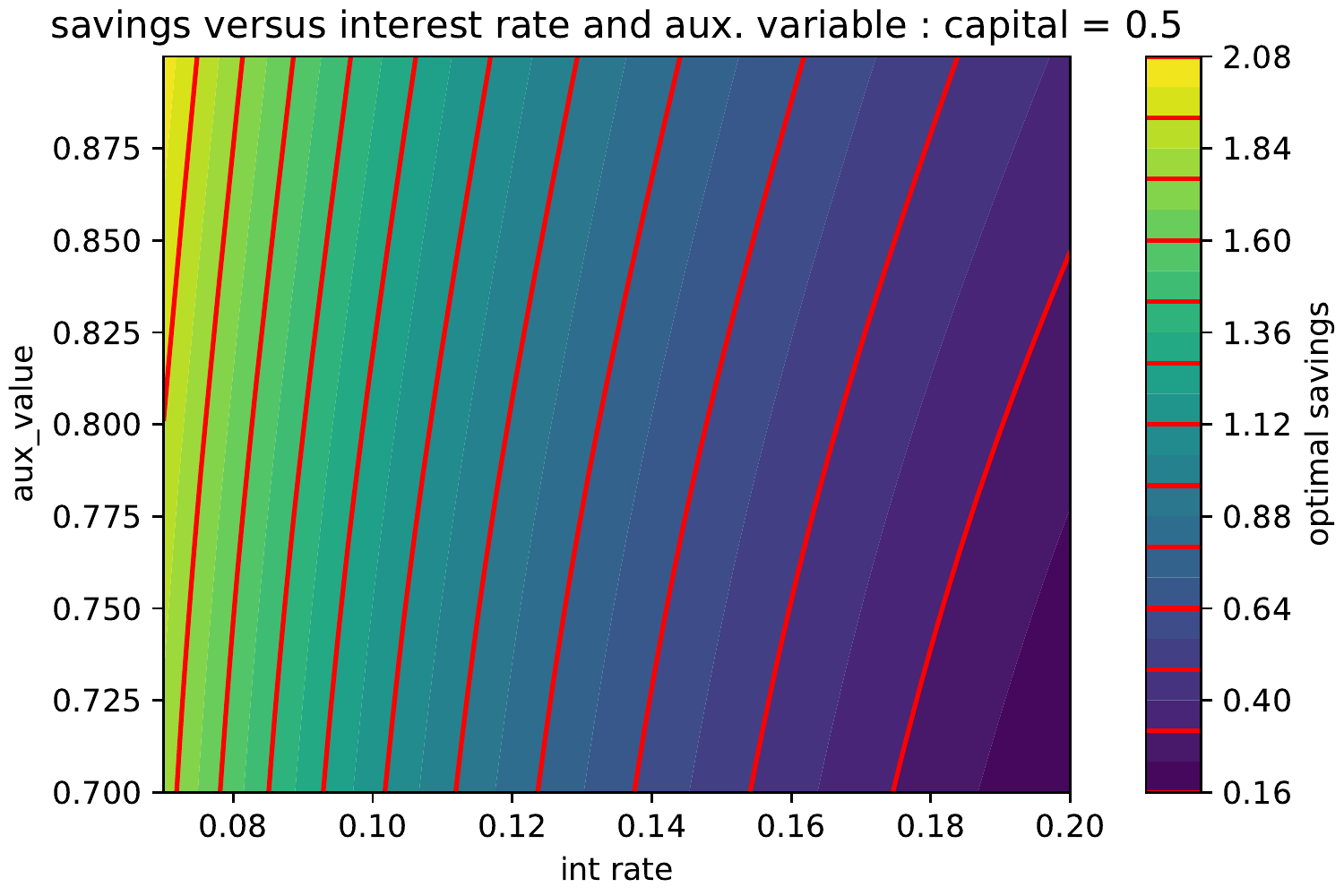}}
   \includegraphics[width=0.4\linewidth]{{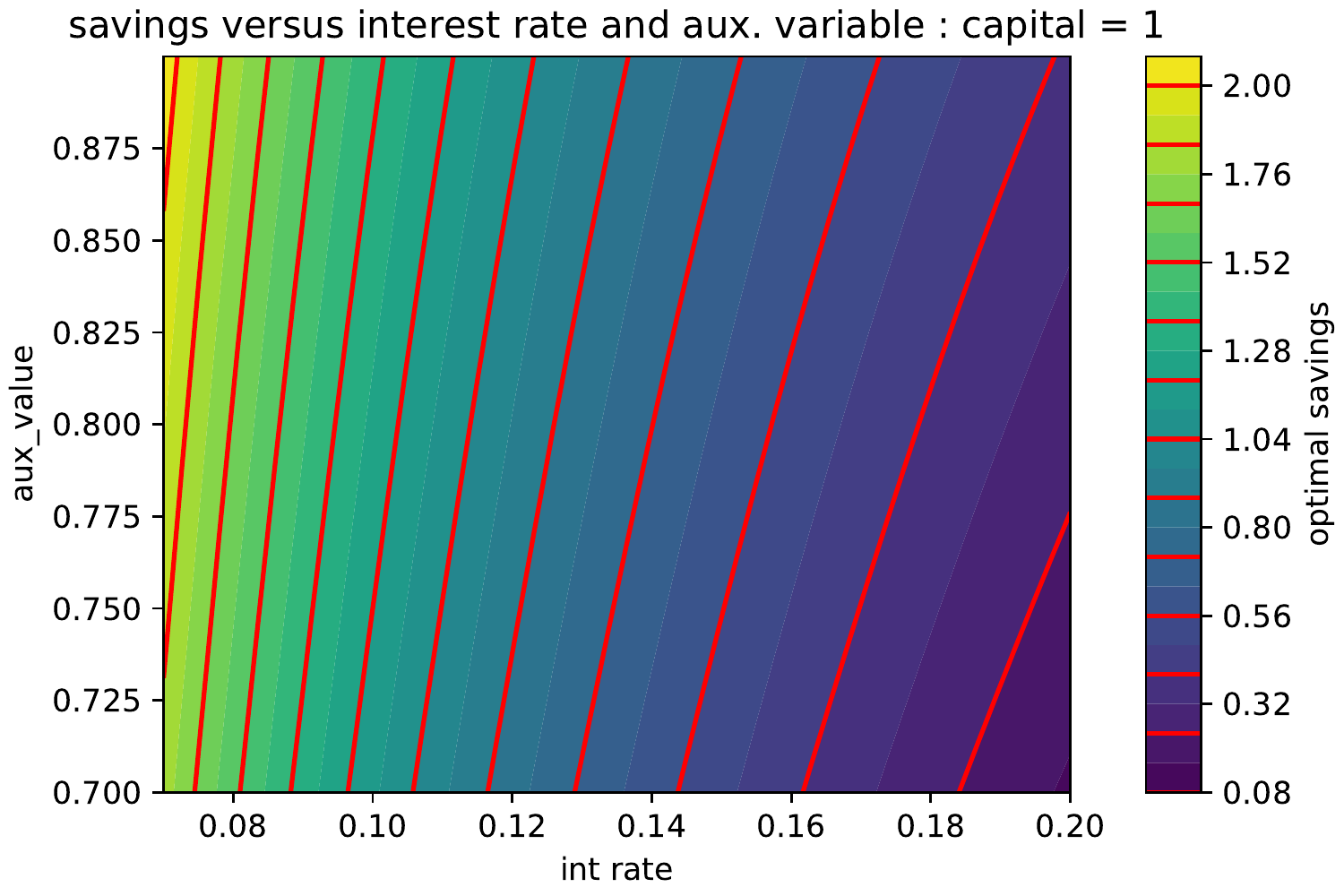}}
 \includegraphics[width=0.4\linewidth]{{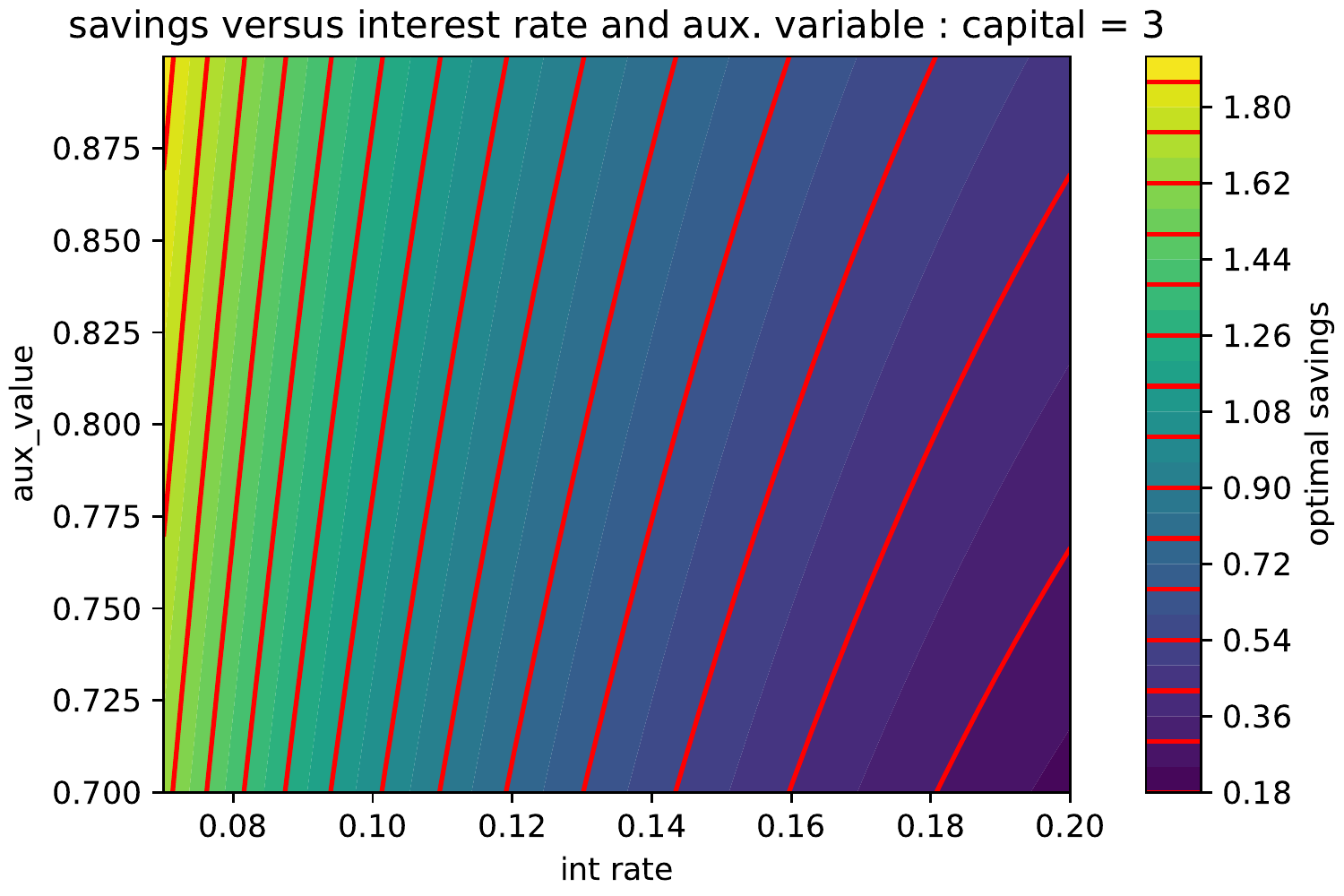}}
 \includegraphics[width=0.4\linewidth]{{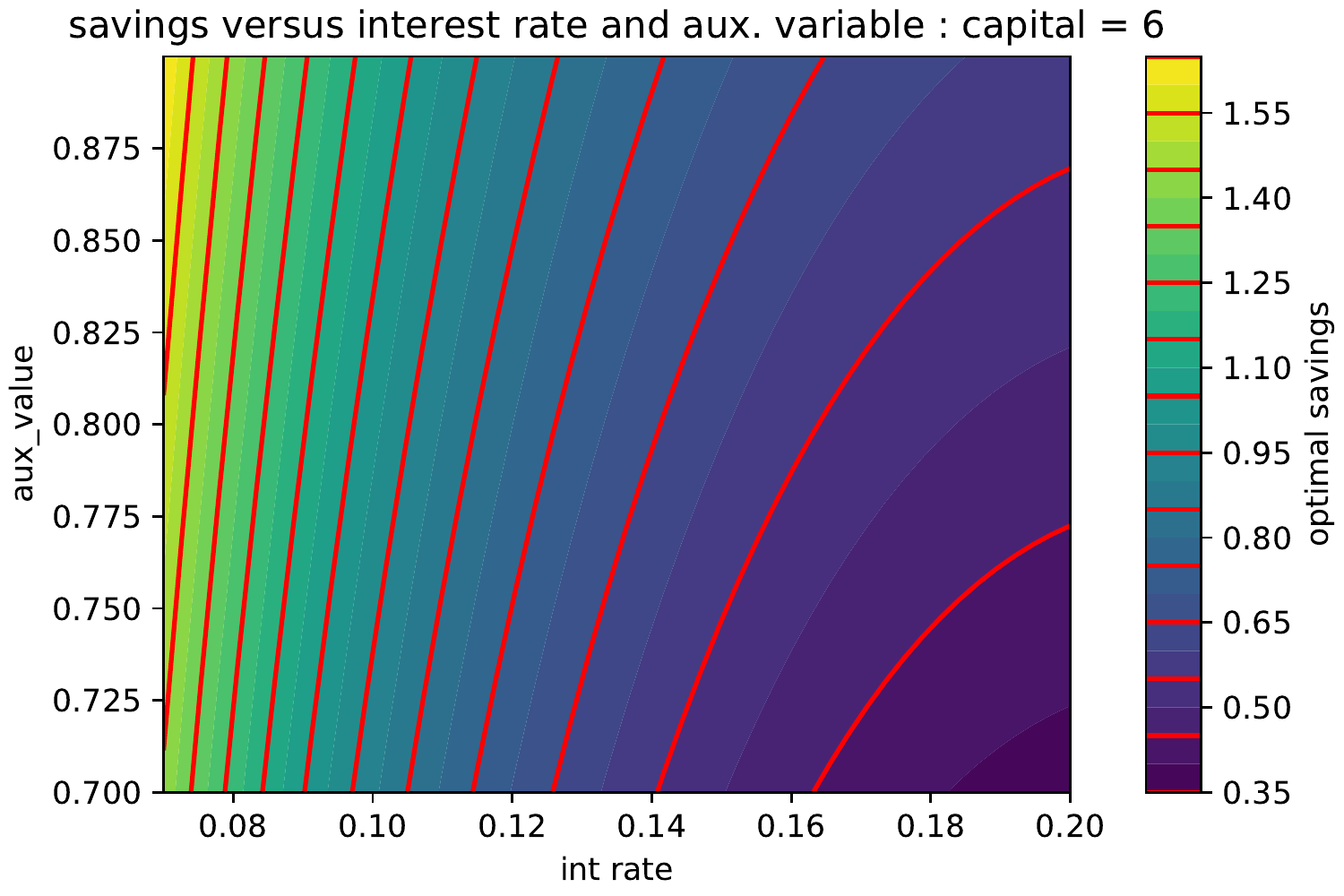}}
 \captionof{figure}{Optimal savings for productive households when $A=A_2$  as a function of  $r$ and the auxiliary variable for $x=0.5,\;1,\;3,\;6$.
  }
  \label{fig:5}
\end{center}

\subsubsection{Correlation between $r$ and the auxiliary variable}

\begin{center}
  \includegraphics[width=0.4\linewidth]{{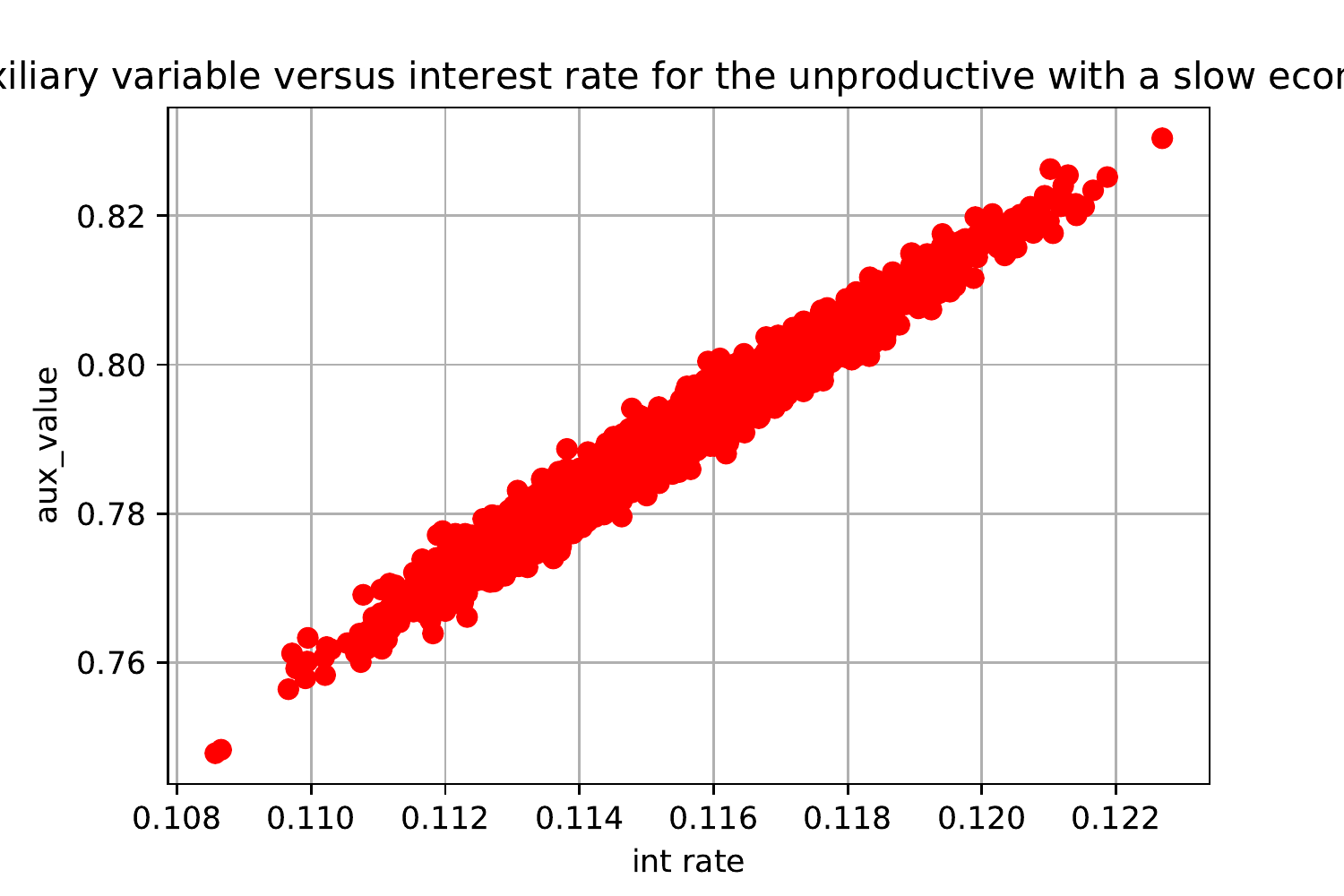}}
  \includegraphics[width=0.4\linewidth]{{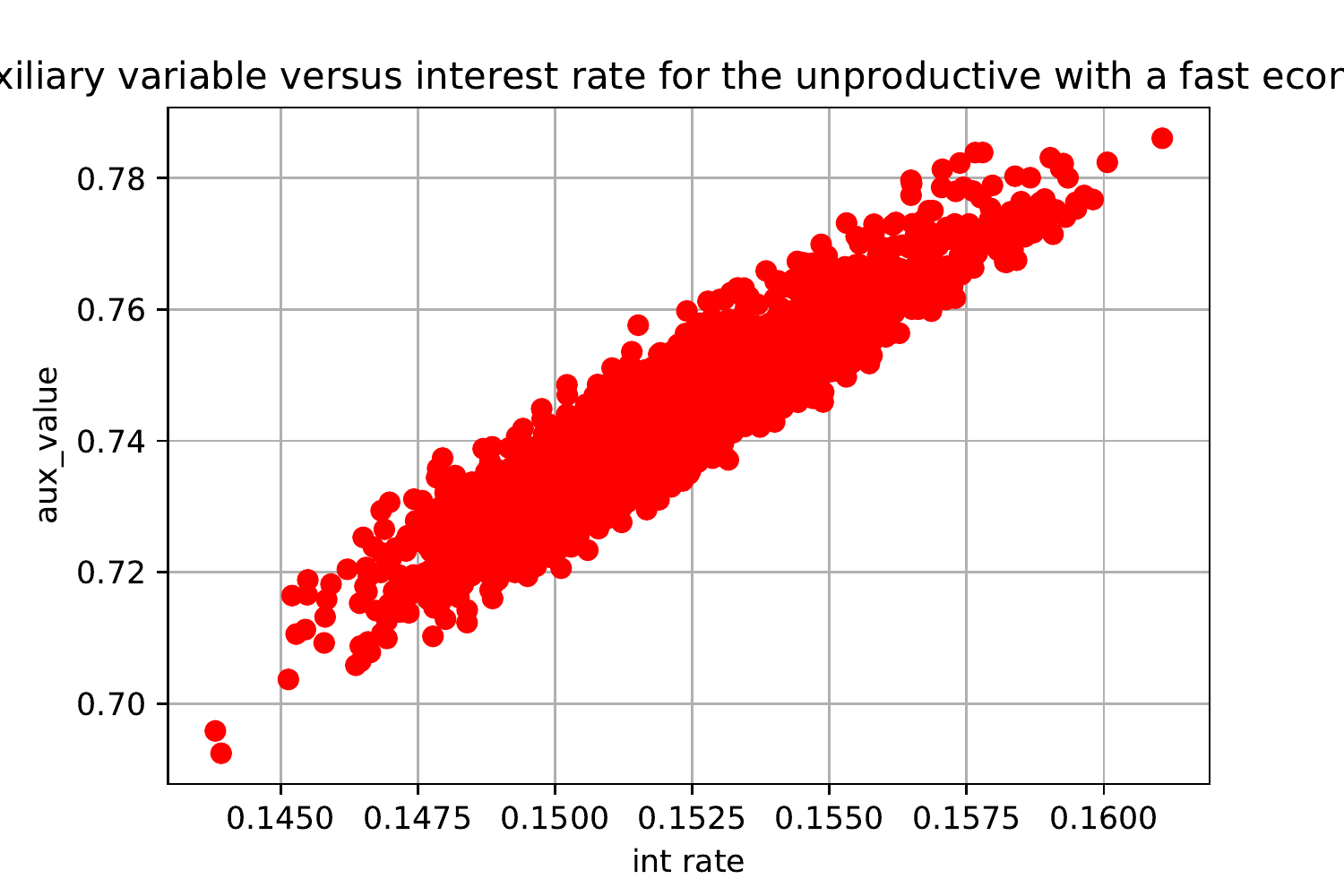}}
  \includegraphics[width=0.4\linewidth]{{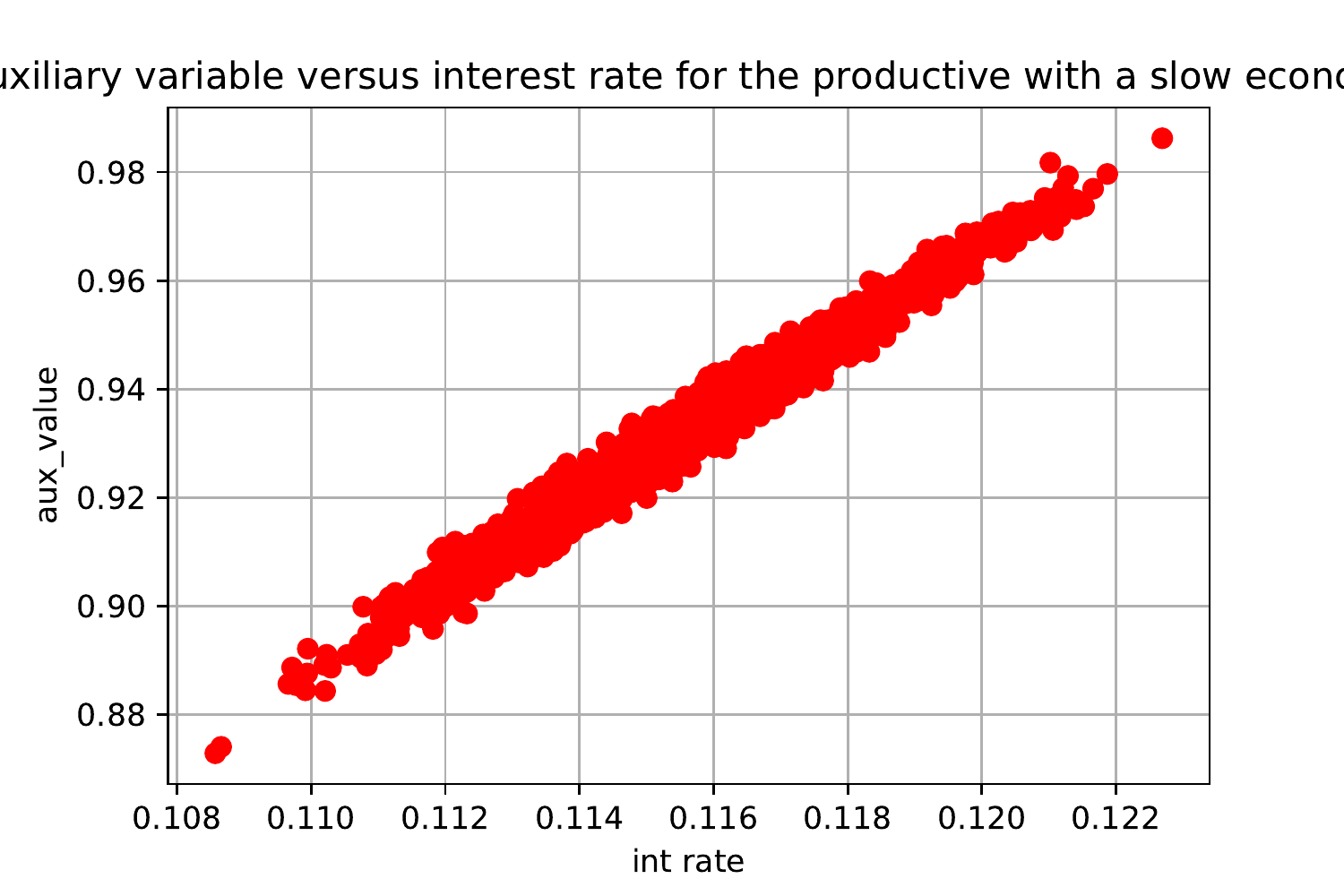}}
   \includegraphics[width=0.4\linewidth]{{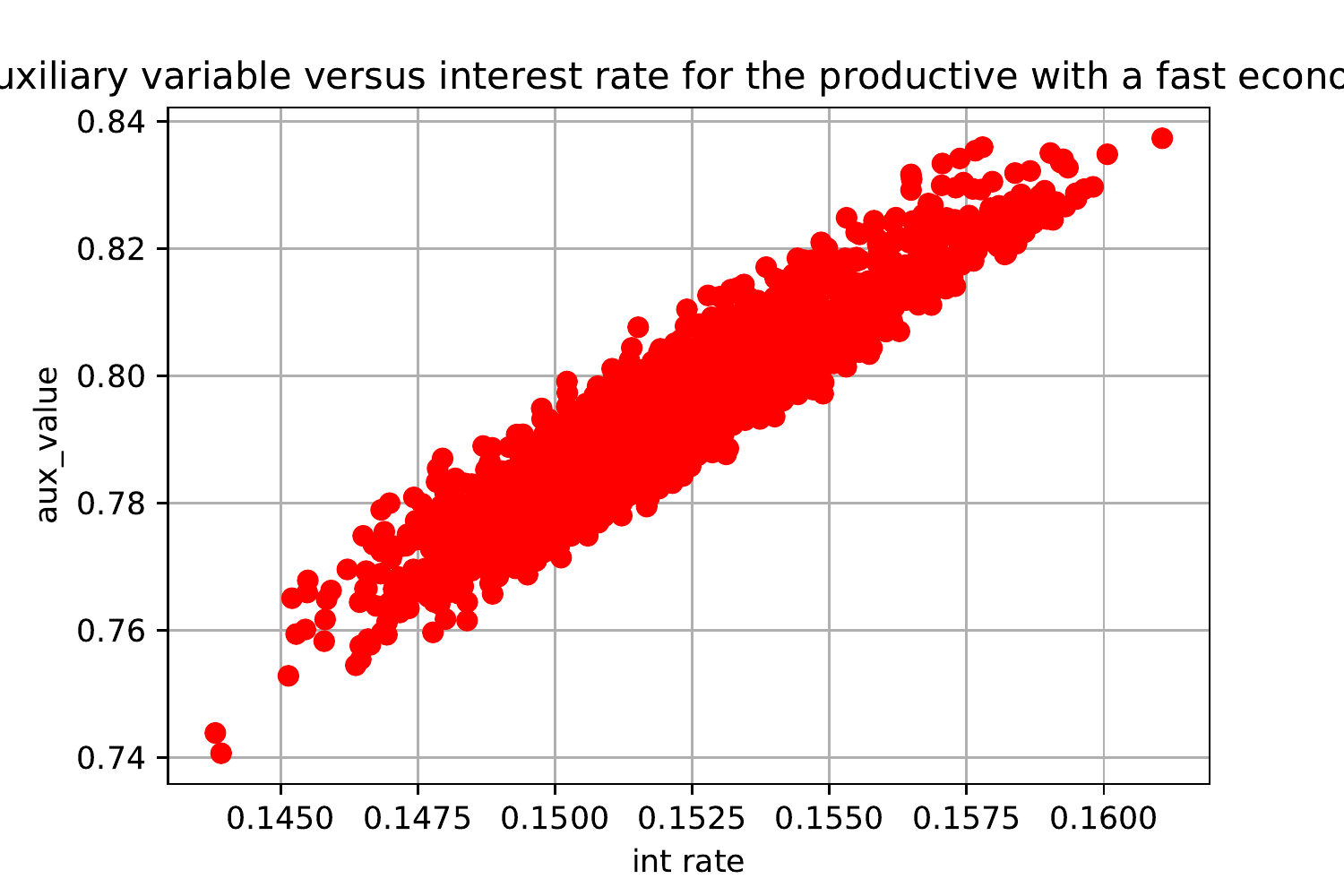}}

  \captionof{figure}{ $\cF_1(m)$ versus $r$ for a sample of measures $m$. Top Left: unproductive and $A=A_1$. Top Right: unproductive and $A=A_2$. Bottom Left: productive and $A=A_1$. Bottom Right: productive and $A=A_2$. We see a rather strong correlation between $r$ and $\cF_1(m)$.
}
  \label{fig:6}
\end{center}

\subsubsection{A tentative conclusion}
In the simulations reported above,  in each situation, ($i=1,2$, $j=1,2$), the approximate solution of the master equation depends on $m$ through the interest rate and an additional variable $\cF_{1,i,j}(m)$,
(the function $\cF_{1,i,j}$ is obtained as a sublayer contained in the first layer of the neural network). While the savings strategies of the households 
mostly depend on the interest rate $r$ as conjectured by Krusell and Smith,  the additional variables $\cF_{1,i,j}(m)$ seem to bring a  significant correction.

\bibliographystyle{siamplain}
\bibliography{KS}

  \end{document}